\documentclass{amsart}

\usepackage{amsmath, amssymb, amsthm}
\usepackage{array}
\usepackage{placeins}
\usepackage[colorlinks=true, linkcolor=blue, citecolor=blue,
            urlcolor=blue]{hyperref}

\theoremstyle{plain}
\newtheorem{theorem}{Theorem}[section]
\newtheorem{proposition}[theorem]{Proposition}
\newtheorem{lemma}[theorem]{Lemma}
\newtheorem{corollary}[theorem]{Corollary}
\newtheorem{conjecture}[theorem]{Conjecture}

\theoremstyle{definition}
\newtheorem{definition}[theorem]{Definition}
\newtheorem{example}[theorem]{Example}

\theoremstyle{remark}
\newtheorem{remark}[theorem]{Remark}

\DeclareMathOperator{\Sing}{Sing}
\newcommand{\PP}{\mathbb{P}}

\renewcommand{\AA}{\mathbb{A}}
\newcommand{\cO}{\mathcal{O}}
\newcommand{\fm}{\mathfrak{m}}
\newcommand{\restr}[2]{\left.{#1}\right|_{#2}}
\newcommand{\crit}{\operatorname{Crit}}
\newcommand{\branch}{\operatorname{Branch}}
\newcommand{\PC}{\operatorname{PC}}
\newcommand{\Irr}{\operatorname{Irr}}
\newcommand{\Spec}{\operatorname{Spec}}

\usepackage{xcolor}

\newcommand{\skill}[1]{{\ttfamily\hyphenchar\font=`\-\relax #1}}

\newcommand{\stacks}[2][]{%
  \cite[\href{https://stacks.math.columbia.edu/tag/#2}{Tag~#2}#1]{Stacks}}

\makeatletter
\let\pcforigbib\thebibliography
\renewcommand{\thebibliography}[1]{\pcforigbib{#1}\raggedright}
\makeatother

\usepackage{etoolbox}
\makeatletter
\let\pcforigbibitem\bibitem
\def\pcfbibitem[#1]#2{%
  \ifstrequal{#2}{Lean4}{\pcforigbibitem[lean]{#2}}{%
    \ifstrequal{#2}{Stacks}{\pcforigbibitem[stacks]{#2}}{%
      \ifstrequal{#2}{MathlibRev26}{\pcforigbibitem[mathlib]{#2}}{%
        \pcforigbibitem[#1]{#2}}}}}
\def\bibitem{\@ifnextchar[\pcfbibitem\pcforigbibitem}
\makeatother

\title{PCF all the way down}

\author{Taylor Dupuy}
\address{Department of Mathematics and Statistics, University of Vermont}
\author{Alon Levy}
\author{Thomas Tucker}
\author{Zhuoran Wang}
\address{Department of Mathematics, University of Rochester}

\date{\today}

\begin{document}

\begin{abstract}
%% >>> [agent-codex 2026-10-26 char0-divisors]
We prove that a finite surjective post-critically finite self-map of a
variety in characteristic zero remains post-critically finite on every
preperiodic integral effective divisor after passing to its periodic
image with its return map.
Iterating this result proves PCF all the way down.

A detailed report on AI usage is contained in the manuscript.
%% <<< [agent-codex 2026-10-26 char0-divisors]
\end{abstract}

\maketitle

	\setcounter{tocdepth}{1}
	\tableofcontents

%% -------------------------------------------------------
\section{Introduction}
\label{sec:intro}
%% -------------------------------------------------------

Let $f \colon X \to X$ be a finite surjective morphism of a variety over a
field of characteristic zero. Its \emph{branch locus} $\branch(f)$ is the
union of the integral effective divisors of $X$ over which $f$ ramifies, and $f$
is
\emph{post-critically finite} (PCF) when the forward orbit
$\bigcup_{n \geq 0} f^{n}(\branch(f))$ consists of finitely many integral
effective
divisors. For a rational map of $\PP^{1}$ this is the classical condition
that every critical point has finite forward orbit.

If $B$ is an integral effective divisor of $X$ with $f^{s}(B) = B$, then $f^{s}$
restricts to a self-map of $B$, and one can ask whether post-critical
finiteness descends to it. Repeating the question for the periodic integral
effective
divisors of $B$, and so on through the dimensions, one arrives at the
question of whether PCF maps are \emph{PCF all the way down}. The question
was posed at the
\href{https://aimath.org/pastworkshops/finitedynamics.html}{AIM workshop
\emph{Postcritically finite maps in complex and arithmetic dynamics}}
held at the American Institute of Mathematics in March~2014
\cite[workshop held March 3--7, 2014]{AIM14}.

Our main result answers the question in characteristic zero.

\begin{theorem}[= Theorem~\ref{thm:general}]
\label{thm:A}
% >>> [agent-claude 2026-10-29 lean-xref]
% Lean: PaperCharZero.theoremA
% <<< [agent-claude 2026-10-29 lean-xref]
Let $X$ be a variety over a field of characteristic zero, let
$f \colon X \to X$ be a finite surjective post-critically finite
morphism, let $A$ be an integral effective divisor of $X$, and let $r \geq 0$
and $s \geq 1$
be integers
with $f^{r+s}(A) = f^{r}(A)$. Then $\restr{f^{s}}{f^{r}(A)}$ is a
post-critically finite self-map of $f^{r}(A)$.
\end{theorem}

Since the restricted map is again a finite surjective self-map of a
variety, Theorem~\ref{thm:A} can be applied to it in turn, and an
induction on the dimension shows that every finite surjective PCF map in
characteristic zero is PCF all the way down in the sense of
Definition~\ref{def:pcfawd} (Corollary~\ref{cor:pcfawd}).

%% >>> [agent-claude 2026-10-04 branch-loci]
Here the branch locus is the \emph{normalized} branch locus
(Definition~\ref{def:crit-branch}), which measures the local degree of $f$
on \emph{analytic branches}.

For a finite surjective morphism $f \colon X \to Y$ and an integral
effective divisor $E$ of $Y$, one usually says that $f$ is
\emph{branched} over $E$ when the fibre of $f$ over the generic point of
$E$ is not reduced. In characteristic zero these divisors are precisely the
codimension-one irreducible components of the image in $Y$ of the support of
$\Omega_{X/Y}$, matching the description of unramified morphisms in the
Stacks Project \stacks[(3)]{02GF}.
%% >>> [agent-claude 2026-10-13 consistency]
In this paper we pass to the normalization: $\branch(f)$ is the image
of the corresponding locus of the lift $\widetilde{f}$ of $f$ to
normalizations (Definition~\ref{def:crit-branch}), which has better
dynamical properties. 
% >>> [agent-codex 2026-11-04 public-references]
Example~\ref{ex:pinched} gives a characteristic-zero
counterexample for the unnormalized convention.
% <<< [agent-codex 2026-11-04 public-references]
%% <<< [agent-claude 2026-10-13 consistency]

The proof has two parts. The first, Theorem~\ref{thm:bounding}, is a
bound for the branch locus of the restriction $g = \restr{f}{A}$ of a
finite surjective map $f \colon X \to Y$ of normal varieties to an integral
effective
divisor $A$ with image $B$: away from the other branch components of $f$,
the singular locus of $Y$ and the singular locus of $B$, the map $g$ is
\'etale. Zariski--Nagata purity makes $f$ \'etale off its branch
divisor, and after shrinking, Abhyankar's lemma presents $f$
\'etale-locally as a disjoint union of Kummer covers $z^{e} = ua$ over a
regular divisor $a = 0$; restricted to the divisor, such a cover is an
isomorphism. The second part, carried out in
Sections~\ref{sec:main}--\ref{sec:general}, uses the bound to confine the
orbit of every branch component of the restriction to a single proper
closed subset of $B$ chosen before any iterate is taken, which leaves room
for only finitely many integral effective divisors. Normality of $X$ is removed
by
passing to the normalization, and periodicity in place of invariance by
passing to an iterate.

% >>> [agent-codex 2026-10-27 three-paper-motivation]
% Finite ramification in fields generated by iterated inverse images is
% another motivation for studying PCF maps. The separate companion
% manuscript proves a finite-ramification theorem under its stated
% hypotheses~\cite[Theorem ``Finite ramification of preimage fields'']%
% {DTWPreimage26}.
% >>> [agent-codex 2026-10-28 defer-preimage-discussion]
% Perhaps surprisingly, finite ramification of iterated inverse-image
% fields does not follow merely by applying PCF all the way down: that
% theorem controls return maps on subvarieties, while the arithmetic
% argument must control the entire inverse-image tower. A separate
% companion manuscript proves finite ramification under its stated
% hypotheses~\cite[Theorem ``Finite ramification of preimage fields'']%
% {DTWPreimage26}.
% <<< [agent-codex 2026-10-28 defer-preimage-discussion]
% <<< [agent-codex 2026-10-27 three-paper-motivation]

%% >>> [agent-codex 2026-10-26 char0-divisors]
Characteristic zero is used in two places: ramification at
codimension-one points is tame, so Abhyankar's lemma applies, and residue
extensions there are separable, so ramification index one means \'etale
(Lemma~\ref{lem:dictionary}).
%% <<< [agent-codex 2026-10-26 char0-divisors]

\subsection{Relation to the literature}

% >>> [agent-codex 2026-11-04 citation-locators]
Over $\mathbb{C}$ the phenomenon was studied earlier. Forn\ae ss and
Sibony \cite[pp.~245--260]{FS92} introduced critically finite maps
of $\PP^{2}$,
and Jonsson \cite[Lemma~1.6 and Remark~1.8]{Jon97} proved the surface case
for an invariant curve; \cite[Remark~1.10]{Jon97} already runs the
mechanism of Theorem~\ref{thm:A} in dimension two for a possibly singular
invariant curve, so in dimension two we claim no priority. In arbitrary
dimension Astorg \cite[Definitions~4--5 and Theorem~1]{Ast20} defined
(weakly) PCF all the way down for endomorphisms of $\PP^{k}(\mathbb{C})$
and proved that a PCF endomorphism is weakly PCF all the way down provided
the intersections of the periodic postcritical components are smooth.
Theorem~\ref{thm:A} removes that smoothness hypothesis from the inheritance
of post-critical finiteness, in the algebraic category, for singular
ambient varieties and for non-normal periodic integral effective divisors; it
says nothing about
Astorg's identification of the descent strata with intersections of
postcritical components, which is what his hypothesis buys. Le
\cite[Theorem~1.1]{Le23} studies the eigenvalues of maps that are weakly
PCF all the way down.
% <<< [agent-codex 2026-11-04 citation-locators]

% >>> [agent-codex 2026-10-28 prior-uses-of-purity]
Purity and Abhyankar's lemma have also been used together to study
tame covers (see section~\ref{sec:provenance}). Bell--Booher--Chen--Liu use purity to reduce an
\'etaleness check to codimension-one local rings
\cite[proof of Thm.~2.2.2]{BBCL}. They then use Abhyankar's lemma
along a normal-crossings boundary to deduce tameness of
restrictions to special fibres
\cite[\S2.3, Lem.~2.3.1]{BBCL}. Kerz--Schmidt compare notions of
tameness in higher dimension and explain the role of a
normal-crossings boundary
\cite[Introduction and Thm.~1.1]{KS10}. Our bounding argument
applies these local tools on a regular open set away from the
other branch components, without a global normal-crossings boundary.
% <<< [agent-codex 2026-10-28 prior-uses-of-purity]

Sections~\ref{sec:defs}
and~\ref{sec:pcf} set up ramification and PCF maps,
Section~\ref{sec:examples} works out examples, Section~\ref{sec:bounding}
proves the bounding theorem, and Sections~\ref{sec:main}
and~\ref{sec:general} prove Theorem~\ref{thm:A}.
%% >>> [agent-codex 2026-10-26 char0-divisors]
Example~\ref{ex:pinched} shows
that Theorem~\ref{thm:A} fails for the unnormalized branch loci, even
in characteristic zero.
%% <<< [agent-codex 2026-10-26 char0-divisors]

\subsection{Conventions}

%% >>> [agent-codex 2026-10-26 char0-divisors]
Throughout, $k$ is a field of characteristic zero. A \emph{variety}
is an integral separated scheme of finite type over $k$; we do not
assume $k$ algebraically closed.
%% <<< [agent-codex 2026-10-26 char0-divisors]

For a variety $V$ we write $\nu = \nu_{V} \colon \widetilde{V} \to V$ for its
normalization, which is finite, and $\Sing(V)$ for its non-regular locus,
which is closed \stacks{07R5}. Closed subsets are given their reduced
induced scheme structure, and a reduced divisor is identified with its
support. We write $\Irr(V)$ for the set of irreducible components of a
closed set~$V$.

\subsection{Software and AI Discussion}
This manuscript is completely formalized (\S\ref{sec:formalization}). 

% >>> [agent-codex 2026-10-27 char0-ai-summary]
OpenAI GPT-5.4 in ChatGPT Classic
\cite[\href{https://openai.com/index/gpt-5-4-thinking-system-card/}{system card}]{OpenAI26},
Anthropic Claude Fable~5
\cite[\href{https://www.anthropic.com/news/claude-fable-5-mythos-5}{model page}]{Ant26}
and Opus~5.5
\cite[\href{https://www.anthropic.com/claude-opus-5-5}{model page}]{AntOpus55}
in Claude Code, and OpenAI GPT-6 Astra
\cite[\href{https://developers.openai.com/api/docs/models/gpt-6-astra}{model page}]{OpenAIAstra26}
and a GPT-6-family model
\cite[\href{https://developers.openai.com/api/docs/guides/latest-model}{model-family guide}]{OpenAIGPT6Family26}
in Codex desktop
\cite[\href{https://openai.com/codex/}{harness page}]{OpenAICodex26} assisted with the
bounding argument, counterexamples, source checks, drafting,
and Lean formalization.

% >>> [agent-codex 2026-10-28 prior-uses-of-purity]
% The techniques of Theorem~\ref{thm:bounding} combining Zariski-Nagata
% Purity and Abhyankar's Lemma showed up in many places.
% The sourced comparison now appears under Relation to the literature.
% <<< [agent-codex 2026-10-28 prior-uses-of-purity]

% <<< [agent-codex 2026-10-27 char0-ai-summary]

\subsubsection{AI findings and procedure}\label{sec:findings}
This paper started as a conventional paper and was completed using AI. 
We describe the devlopments with AI.
\begin{itemize}
	\item The first AI-assisted result was a counterexample to a statement about how the branch locus behave with respect to restriction on divisors which is now Example~\ref{ex:cusp}.
	\item In searching for a repair several months later a Fable agent suggested replacing
	the failed equality by the upper bound of
	Theorem~\ref{thm:bounding} and using purity followed by
	Abhyankar's lemma.
	This was enough to complete the proof.
	\item We then searched to see if the bound was actually attained. This produced Example~\ref{ex:branch-singular},
	which shows the ambient-singularity term is needed in the
	bounding theorem.
	% >>> [agent-codex 2026-10-27 pinched-provenance-transfer]
	%% Editorial note addressed in the companion article.
	% The agent provenance for Example~\ref{ex:pinched} now appears in
	% the companion article; the counterexample itself remains here.
	% <<< [agent-codex 2026-10-27 pinched-provenance-transfer]
\end{itemize}
Some other examples were found by AI searches.

The counter-example to a naive claim was done in an interactive session with OpenAI GPT-5.4 in ChatGPT Classic
\cite[\href{https://openai.com/index/gpt-5-4-thinking-system-card/}{GPT-5.4 system card}]{OpenAI26}.
The bounding Lemma was constructed by Claude Fable~5 in Claude Code
\cite[\href{https://www.anthropic.com/news/claude-fable-5-mythos-5}{Fable~5 model page}]{Ant26}.
The example providing some sharpness to our bounding lemma used OpenAI GPT-6 Astra in Codex desktop at effort xhigh.
\cite[\href{https://developers.openai.com/api/docs/models/gpt-6-astra}{Astra model page}]{OpenAIAstra26}
\cite[\href{https://openai.com/codex/}{Codex harness page}]{OpenAICodex26}.

\subsubsection{Provenance}\label{sec:provenance}
% >>> [agent-codex 2026-10-27 first-proof-provenance]
%% \taylor{Cite 1st Proof} -- addressed by the pinpoint citation below.
The authors checked whether the purity--Abhyankar strategy had
antecedents in the literature. This check was motivated in part by
reported instances of AI solutions reusing earlier work without
citation~\cite[p.~5]{ASWW26}. 

A targeted Codex literature search
located related arguments in SGA~1, Expos\'e~X
\cite[Th\'eor\`eme~3.1 and Lemme~3.6]{SGA1},
Grothendieck--Murre \cite[Chapters~1--5]{GM71},
Helminck \cite[Proposition~1]{Hel17},
Kerz--Schmidt \cite[Lemmas~3.5--3.6]{KS10}, and
Carvajal-Rojas--St\"abler \cite[\S1.1 and Theorem~3.7]{CRS23}.
The authors reviewed these items.
Although the original search yielding the bounding lemma did not cite these items, it shows the trick the agent supplied is not new.
% <<< [agent-codex 2026-10-27 first-proof-provenance]

The cone underlying Example~\ref{ex:branch-singular} occurs in
Reid's \cite[Example~1.2]{ReidCyclic}
\href{https://mreid.warwick.ac.uk/surf/more/cyclic.pdf}{notes};
the restriction calculation supplied here appears to be new.

% >>> [agent-codex 2026-11-04 public-references]
Example~\ref{ex:pinched} was developed with Claude Opus~5.5 in
Claude Code and reviewed by Claude Fable~5.1 subagents.
% <<< [agent-codex 2026-11-04 public-references]

\subsubsection{Due diligence}
% >>> [agent-codex 2026-10-28 due-diligence-companions]
Theorem~\ref{thm:general} concerns integral effective divisors in
characteristic zero. 
We explored the natural generalizations so as to provide a natural stopping point for ourselves and a natural starting point for future research (cleaning-up potential ``one-shots'' by modern frontier agents). 

In particular, 
% >>> [agent-codex 2026-11-04 public-references]
\begin{itemize}
	\item We explored two natural generalizations which will appear in a companion article in preparation~\cite{DTW26PCFAWD}.
	The first is taking the base field to be positive characteristic and the second is restriction to
	preperiodic subvarieties of any codimension.
	These could not be achieved in a ``one-shot'' by a frontier agent. 
	This manuscript also discusses replacing the divisorial locus with the non-\'{e}tale
	locus or the standard branch locus.
	\item The natural questions about the rigid analytic and complex
	analytic cases appear to be more or less ``one-shot'' by frontier agents and we will be providing separate AI-generated supplements for these cases: 
	\cite{DLTW26Rigid,DLTW26Complex}.
	These will be released to hexagon.
\end{itemize}
% <<< [agent-codex 2026-11-04 public-references]

% <<< [agent-codex 2026-10-28 due-diligence-companions]

\subsubsection{Other}
Appendix~\ref{app:ai} describes the methods (typesetting, agent management, record keeping, natural language proofs, adversarial reviews, AI assisted computations and this statement), proofreading and adversarial review, and the formalization (\S\ref{sec:formalization}). Appendix~\ref{app:tokens} reports token accounting.

%% -------------------------------------------------------
\section{Ramification, critical loci and branch loci}
\label{sec:defs}
%% -------------------------------------------------------

Let $g \colon A \to B$ be a finite surjective morphism of varieties. It
lifts uniquely to a morphism
$\widetilde{g} \colon \widetilde{A} \to \widetilde{B}$ of normalizations
with $\nu_{B} \circ \widetilde{g} = g \circ \nu_{A}$
\stacks[(4)]{035Q}. The lift is finite, because $g \circ \nu_{A}$ is finite
and $\nu_{B}$ is separated \stacks[(1)]{035D}, and surjective, because its
image is closed and maps onto $B$, so has dimension $\dim \widetilde{B}$.

%% >>> [agent-claude 2026-10-04 taylor-tags]
%% \taylor{make this $g$} -- addressed: the map is now $g$.
%% <<< [agent-claude 2026-10-04 taylor-tags]
%% [agent-claude 2026-10-13 consistency] Definition def:branch-loci and
%% the remark on the etale branch locus moved to Appendix A, where
%% Branch_0 and Branch_1 are now defined; Branch is defined only by
%% Definition def:crit-branch below.

We will connect the definition of the branch and critical locus to ramification
degrees in a moment.

\begin{lemma}[Lifting integral effective divisors]
	\label{lem:lift}
	% >>> [agent-claude 2026-10-29 lean-xref]
	% Lean: Pcf.finiteDominantFiniteType_divisor_lift,
	%   Pcf.finiteDominantFiniteType_divisor_lift_into
	% <<< [agent-claude 2026-10-29 lean-xref]
	Let $\pi \colon V \to W$ be a finite surjective morphism of varieties and
	let $Q$ be an integral effective divisor of $W$. Then there is an integral
	effective divisor $P$ of $V$ with
	$\pi(P) = Q$, and only finitely many integral effective divisors of $V$ map
	into $Q$.
\end{lemma}

\begin{proof}
	The closed set $\pi^{-1}(Q)$ has finitely many irreducible components,
	whose images are closed and cover the irreducible set $Q$; so one
	component $P$ maps onto $Q$, and $\dim P = \dim Q$ because $\pi$ is
	finite, whence $P$ is an integral effective divisor as in
	Lemma~\ref{lem:restriction}. An integral effective divisor
	of $V$ mapping into $Q$ lies in $\pi^{-1}(Q)$ and has dimension
	$\dim Q$, which is the largest dimension of a component of
	$\pi^{-1}(Q)$, so it is one of those components.
\end{proof}

%% >>> [agent-claude 2026-10-04 taylor-tags]
%% \taylor{track down pinpoint Atiyah-MacDonald citation} -- addressed:
%% cited to the Stacks Project definition, which was checked; no
%% Atiyah--Macdonald pinpoint could be verified.
%% \taylor{Add tildes to all the normal schemes} -- addressed.
We use the ramification index of an extension of discrete valuation rings
\stacks{09E4}. Let $\widetilde{g} \colon \widetilde{A} \to \widetilde{B}$
be a finite surjective morphism of normal varieties, let $\widetilde{D}$
be an integral effective divisor of $\widetilde{A}$, and let
$\widetilde{E} = \widetilde{g}(\widetilde{D})$, an integral effective
divisor of $\widetilde{B}$ (Lemma~\ref{lem:restriction}). The local rings
$\cO_{\widetilde{B},\widetilde{E}}$ and $\cO_{\widetilde{A},\widetilde{D}}$
at the generic points are discrete valuation rings, being normal of
dimension one, and $\cO_{\widetilde{B},\widetilde{E}} \to
\cO_{\widetilde{A},\widetilde{D}}$ is an extension of discrete valuation
rings; we write $e(\widetilde{D}/\widetilde{E})$ for its ramification
index. We say that $\widetilde{g}$ is \emph{ramified} at $\widetilde{D}$,
or that $\widetilde{D}$ is \emph{critical}, and that $\widetilde{E}$ is
\emph{branched}, if $e(\widetilde{D}/\widetilde{E}) > 1$. As the fibre of
$\widetilde{g}$ over the generic point of $\widetilde{E}$ is
$\prod_{\widetilde{D}}
\cO_{\widetilde{A},\widetilde{D}}/\fm_{\widetilde{D}}^{e(\widetilde{D}/\widetilde{E})}$,
the product over the integral effective divisors $\widetilde{D}$ over
$\widetilde{E}$, %% >>> [agent-claude 2026-10-13 consistency]
the fibre over the generic point of $\widetilde{E}$ is reduced exactly
when $\widetilde{E}$ is not branched.
%% <<< [agent-claude 2026-10-13 consistency]
%% <<< [agent-claude 2026-10-04 taylor-tags]

\begin{definition}[Normalized critical and branch loci]
	\label{def:crit-branch}
	% >>> [agent-claude 2026-10-29 lean-xref]
	% Lean: Pcf.PaperCharZero.schemeNormalizedCriticalDivisors,
	%   Pcf.schemeNormalizedBranch
	% <<< [agent-claude 2026-10-29 lean-xref]
	Let $g \colon A \to B$ be a finite surjective morphism of varieties.
The \emph{normalized critical locus}, or simply the
		\emph{critical locus}, $\crit(g) \subseteq A$ is the union of
		the integral effective divisors $\nu_{A}(\widetilde{D})$ for
		$\widetilde{D}$ critical for $\widetilde{g}$.
\end{definition}

The \emph{normalized branch locus}, or simply the \emph{branch
	locus}, is $\branch(g) = g(\crit(g)) \subseteq B$.

By construction $\crit(g) = \nu_{A}(\crit(\widetilde{g}))$ and
\begin{equation}\label{eq:branch-normalization}
% >>> [agent-claude 2026-10-29 lean-xref]
% Lean: Pcf.schemeNormalizedBranch_eq_normalization_image
% <<< [agent-claude 2026-10-29 lean-xref]
	\branch(g) = g(\nu_{A}(\crit(\widetilde{g})))
	= \nu_{B}(\widetilde{g}(\crit(\widetilde{g})))
	= \nu_{B}(\branch(\widetilde{g})) .
\end{equation}
The critical locus of $\widetilde{g}$ is a finite union of integral effective
divisors,
because an integral effective divisor with $e > 1$ lies in the non-\'etale locus
\stacks[(6)]{02GU}, which is closed, and proper since $\widetilde{g}$ is
\'etale at the generic point; such a divisor is one of its finitely many
components of codimension one.
Hence $\crit(g)$ and $\branch(g)$ are reduced divisors, and their
irreducible components are integral effective divisors by
Lemma~\ref{lem:restriction}.
%% >>> [agent-claude 2026-10-13 consistency]
Equivalently, by the description of the fibres of $\widetilde{g}$ above,
$\branch(g)$ is the image under $\nu_{B}$ of the union of the integral
effective divisors of $\widetilde{B}$ over whose generic point the fibre
of $\widetilde{g}$ is not reduced.

% >>> [agent-codex 2026-10-27 char0-branch-warning]
We use the normalized branch locus because it obeys the exact tower
formula needed below (Remark~\ref{rem:why-normalized}).
Example~\ref{ex:pinched} shows that descent fails for the
nonreduced-fibre and non-\'etale loci on a nonnormal variety.
% <<< [agent-codex 2026-10-27 char0-branch-warning]
%% <<< [agent-claude 2026-10-13 consistency]

%% >>> [agent-claude 2026-10-04 branch-loci]
\begin{remark}[Why normalize]
\label{rem:why-normalized}
% >>> [agent-claude 2026-10-29 lean-xref]
% Lean: Pcf.Examples.CharZeroPaperCuspHomeomorph.normalizationMap_isHomeomorph,
%   Pcf.Examples.CharZeroPaperCuspHomeomorph.normalizationMap_branch_empty
% <<< [agent-claude 2026-10-29 lean-xref]
%% \taylor{Use $D$ and $E$ here} -- addressed.
The normalization is built into the definition because ramification is
attached to analytic branches and not to integral effective divisors of
$A$. What the definition measures is the local degree of $g$ along each
branch: over $\mathbb{C}$, if $\widetilde{D}$ is an integral effective
divisor of $\widetilde{A}$ with $D = \nu_{A}(\widetilde{D})$ and
$\widetilde{E} = \widetilde{g}(\widetilde{D})$, then at a general point of
$D$ the map $g$ has local degree $e(\widetilde{D}/\widetilde{E})$ on every
analytic branch of $A$ coming from $\widetilde{D}$. This is the quantity that a
critical
point records in dynamics, and it does not see singularities: the
normalization of a cuspidal curve is a homeomorphism, and it is not
branched.

% >>> [agent-codex 2026-10-27 char0-branch-warning]
The exact tower formula of Lemma~\ref{lem:tower} lets
post-critical finiteness pass between a map and its normalization
(Lemmas~\ref{lem:normalization} and~\ref{lem:covers}).
The unnormalized alternatives require a different analysis;
Example~\ref{ex:pinched} shows why they cannot replace this one.
% <<< [agent-codex 2026-10-27 char0-branch-warning]
\end{remark}
%% <<< [agent-claude 2026-10-04 branch-loci]

\begin{lemma}[Ramification and \'etaleness]
\label{lem:dictionary}
% >>> [agent-claude 2026-10-29 lean-xref]
% Lean: Pcf.PaperCharZero.ramificationIdx_eq_one_iff_etale_at_divisor_charZero
% <<< [agent-claude 2026-10-29 lean-xref]
Let $g \colon A \to B$ be a finite surjective morphism of normal varieties
over a field of characteristic zero and let $P$ be an integral effective divisor
of $A$ with
generic point $\xi$. Then $e(P/g(P)) = 1$ if and only if $g$ is \'etale
at $\xi$.
\end{lemma}

\begin{proof}
Put $R = \cO_{B,g(\xi)}$ and $S = \cO_{A,\xi}$, both discrete valuation
rings. The ring $S$ is a domain containing $R$, hence a torsion-free and
therefore flat $R$-module \stacks[(1)]{0AUW}. By the equivalence of
parts~(1) and~(6) of \stacks{02GU}, $g$ is \'etale at $\xi$ if and only if
$\fm_{R}S = \fm_{S}$ and $\kappa(\xi)/\kappa(g(\xi))$ is finite separable.
The residue extension is finite because $g$ is finite, and separable
because the characteristic is zero, so the condition is
$e(P/g(P)) = 1$.
\end{proof}

% >>> [agent-codex 2026-10-27 char0-local-index]
\begin{lemma}[\'Etale restriction and normalized ramification]
\label{lem:etale-normalized-index}
% >>> [agent-claude 2026-10-29 lean-xref]
% Lean: Pcf.PaperCharZero.etale_normalized_index
% <<< [agent-claude 2026-10-29 lean-xref]
Let $g\colon A\to B$ be a finite surjective morphism of varieties
over a field of characteristic zero. Let $P$ be an integral effective
divisor of $A$, and suppose that $B$ is regular at the generic point
of $g(P)$ and that $g$ is \'etale at the generic point of $P$.
Then the normalized ramification index along every branch over $P$
is one.
\end{lemma}

\begin{proof}
At the generic point of $g(P)$, the local ring of $B$ is regular.
Since $g$ is \'etale at the generic point of $P$, the local ring
of $A$ there is regular as well. The normalizations of $A$ and
$B$ are therefore isomorphisms near these points. Their lifted
map agrees locally with $g$, so Lemma~\ref{lem:dictionary} gives
ramification index one. There is only one normalization branch
at the generic point of $P$.
\end{proof}
% <<< [agent-codex 2026-10-27 char0-local-index]

%% -------------------------------------------------------
\section{Post-critically finite morphisms}
\label{sec:pcf}
%% -------------------------------------------------------

%% >>> [agent-codex 2026-10-26 char0-divisors]
Let $f \colon X \to X$ be a finite surjective morphism of a variety.
An integral effective divisor $A$ of $X$ is \emph{preperiodic} if
$f^{r+s}(A) = f^{r}(A)$ for some $r \geq 0$ and $s \geq 1$.
It is \emph{invariant} if $f(A) = A$.
%% <<< [agent-codex 2026-10-26 char0-divisors]

\begin{definition}[Postcritical set and PCF]
\label{def:pcf}
% >>> [agent-claude 2026-10-29 lean-xref]
% Lean: Pcf.schemeNormalizedPostcriticalSet, Pcf.IsNormalizedPCFOfFiniteType,
%   Pcf.PaperCharZero.postcriticalSubscheme
% <<< [agent-claude 2026-10-29 lean-xref]
Let $f \colon X \to X$ be a finite surjective morphism of a variety.
\begin{enumerate}
\item The \emph{postcritical set} of $f$ is the set of integral effective
divisors
\[
\PC(f) = \{\, f^{n}(D) : D \in \Irr(\branch(f)),\ n \geq 0 \,\} .
\]
\item The map $f$ is \emph{post-critically finite} (PCF) if $\PC(f)$ is
finite.
\end{enumerate}
\end{definition}

Each $f^{n}(D)$ is an integral effective divisor by Lemma~\ref{lem:restriction},
and the union
of the members of $\PC(f)$ is
$\bigcup_{n \geq 1} f^{n}(\crit(f)) = \bigcup_{n \geq 0}
f^{n}(\branch(f))$, the orbit of the critical locus. When $f$ is PCF this
union is a reduced divisor, which we also denote $\PC(f)$. By definition
$f$ carries each member of $\PC(f)$ to a member of $\PC(f)$.

%% >>> [agent-claude 2026-10-04 branch-loci]

%% >>> [agent-codex 2026-10-26 char0-divisors]
% >>> [agent-codex 2026-10-28 statement-notation-charzero]
\begin{lemma}[Restriction to an integral effective divisor]
	\label{lem:restriction}
	% >>> [agent-claude 2026-10-29 lean-xref]
	% Lean: Pcf.PaperCharZero.integral_effective_divisor_restriction,
	%   Pcf.PaperCharZero.DimensionDrop.divisorClosure_dimension_add_one_eq
	% <<< [agent-claude 2026-10-29 lean-xref]
	Let $f \colon X \to Y$ be a finite surjective morphism of varieties
	and let $A$ be an integral effective divisor of $X$. Then the
	image $f(A)$, with its reduced induced structure, is an integral
	effective divisor of $Y$, and the restriction of $f$ to $A$ is a
	finite surjective morphism onto that image.
\end{lemma}
% <<< [agent-codex 2026-10-28 statement-notation-charzero]
%% <<< [agent-codex 2026-10-26 char0-divisors]

\begin{proof}
	% >>> [agent-codex 2026-10-28 statement-notation-charzero]
	Let $B$ be the reduced closed subscheme supported on the closure
	of $f(A)$ in $Y$.
	Write $\iota_{A}$ and $\iota_{B}$ for the closed immersions. A closed
	immersion is finite \stacks{035C}, so $f \circ \iota_{A}$ is finite, hence
	closed; thus the support of $B$ is $f(A)$ and is irreducible.
	As $A$ is reduced,
	$f \circ \iota_{A}$ factors through $\iota_{B}$ \stacks{0356},
	giving $g = \restr{f}{A} \colon A \to B$;
	% <<< [agent-codex 2026-10-28 statement-notation-charzero]
	it is finite because $f \circ \iota_{A}$ is finite and $\iota_{B}$ is
	separated \stacks[(1)]{035D}, and surjective by construction.

	Since $f$ is surjective, $\dim X = \dim Y$ and $\dim A = \dim B$,
	since finite morphisms are integral and integral surjective morphisms
	preserve dimension \stacks{0ECG}. For an irreducible closed subset $W$ of
	a variety $V$ one has $\dim V = \dim W + \operatorname{codim}(W, V)$
	\stacks[(3), (6), (10)]{0A21}. So $\operatorname{codim}(B, Y) = \dim
	Y - \dim B
	= \dim X - \dim A = 1$.
\end{proof}

We can now state the property of the title. The recursion is on the
dimension, which drops by one at each step by Lemma~\ref{lem:restriction}.

\begin{definition}[PCF all the way down]
	\label{def:pcfawd}
	% >>> [agent-claude 2026-10-29 lean-xref]
	% Lean: Pcf.SchemePCFAllTheWayDown,
	%   Pcf.PaperCharZero.pcfAllTheWayDown_iff_at_dimension_bound,
	%   Pcf.PaperCharZero.dimension_zero_isPCFAllTheWayDown
	% <<< [agent-claude 2026-10-29 lean-xref]
	Let $f \colon X \to X$ be a finite surjective morphism of a variety. If
	$\dim X = 0$, then $f$ is \emph{PCF all the way down}. If $\dim X \geq 1$,
	then $f$ is \emph{PCF all the way down} if $f$ is PCF and, for every
	preperiodic integral effective divisor $A$ of $X$ with $f^{r+s}(A) = f^{r}(A)$,
	$r \geq 0$,
	$s \geq 1$, the map $\restr{f^{s}}{B} \colon B \to B$ with
	$B = f^{r}(A)$ is PCF all the way down.
\end{definition}

In dimension one, PCF all the way down is the same as PCF. The question
of the introduction is the following.

\begin{conjecture}
	\label{conj:pcfawd}
	% >>> [agent-claude 2026-10-29 lean-xref]
	% Lean: PaperCharZero.cor_pcfawd
	% <<< [agent-claude 2026-10-29 lean-xref]
	Every finite surjective post-critically finite self-map of a variety over
	a field of characteristic zero is PCF all the way down.
\end{conjecture}

Corollary~\ref{cor:pcfawd} proves Conjecture~\ref{conj:pcfawd}.

%% <<< [agent-claude 2026-10-04 branch-loci]

%% -------------------------------------------------------
\section{Examples}
\label{sec:examples}
%% -------------------------------------------------------

% >>> [agent-codex 2026-10-27 char0-examples]
The first example shows the recursion on coordinate hyperplanes. The
next two give projective-plane maps with less trivial critical orbits.

\begin{example}[Coordinatewise squaring]
\label{ex:squaring-PN}
% >>> [agent-claude 2026-10-29 lean-xref]
% Lean: Pcf.Examples.CharZeroPaperSquaringBranch.squaringMap_branch,
%   Pcf.Examples.CharZeroPaperSquaringIndex.squaringMap_ramificationIdx_hyperplane,
%   Pcf.Examples.CharZeroPaperSquaringRestriction.hyperplane_restriction_squaring_isPCF,
%   Pcf.Examples.CharZeroPaperSquaringBranch.squaringMap_isPCF,
%   Pcf.Examples.CharZeroPaperProjectiveAllTheWayDown.squaringMap_pcfAllTheWayDown,
%   Pcf.Examples.CharZeroPaperSquaringBranch.diagonal_not_mem_branch_pow,
%   Pcf.Examples.CharZeroPaperSquaringChartJacobian.squaring_chart_jacobian
% <<< [agent-claude 2026-10-29 lean-xref]
For $N\geq2$, let
\[
\psi([x_0:\cdots:x_N])=[x_0^2:\cdots:x_N^2]
\quad\text{on }\PP^N,
\qquad H_i=V(x_i).
\]
On each standard affine chart the Jacobian is a nonzero scalar times
the product of the remaining affine coordinates. Thus the critical
and branch divisors are the $H_i$, each of which is invariant.
The restriction to $H_i\cong\PP^{N-1}$ is again coordinatewise
squaring. The map $\psi$ is PCF, and Corollary~\ref{cor:pcfawd}
shows that it is PCF all the way down, including on periodic
divisors other than the $H_i$.

The integral effective divisor $D=V(x_0-x_1)$ is invariant, since
$x_0^2-x_1^2=(x_0-x_1)(x_0+x_1)$, but $D$ is not a component of any
$\branch(\psi^m)$: those components are coordinate hyperplanes.
Thus periodic integral effective divisors need not be postcritical.
\end{example}

\begin{example}[Belk--Koch map]
\label{ex:koch}
% >>> [agent-claude 2026-10-29 lean-xref]
% Lean: Pcf.Examples.CharZeroPaperPlaneBranch.kochMap_branch,
%   Pcf.Examples.CharZeroPaperPlaneBranch.kochMap_postcritical,
%   Pcf.Examples.CharZeroPaperPlaneBranch.kochMap_isPCF
% <<< [agent-claude 2026-10-29 lean-xref]
Belk and Koch give the PCF endomorphism
\[
F([x:y:z])=[x^2-y^2:x^2-z^2:x^2]
\quad\text{of }\PP^2.
\]
Its postcritical locus is a union of six lines
\cite[Introduction, pp.~1--2]{BelkKoch10}. This supplies a compact
example whose critical orbit differs from the coordinatewise map.
\end{example}

\begin{example}[Dupont's map]
\label{ex:dupont}
% >>> [agent-claude 2026-10-29 lean-xref]
% Lean: Pcf.Examples.CharZeroPaperPlaneBranch.dupontMap_criticalDivisors,
%   Pcf.Examples.CharZeroPaperPlaneBranch.dupontMap_postcritical,
%   Pcf.Examples.CharZeroPaperPlaneBranch.dupontMap_isPCF,
%   Pcf.Examples.CharZeroPaperPlaneMinor.dupontLinearMatrix_det,
%   Pcf.Examples.CharZeroPaperPlaneMinor.dupontMap_H_isStrictlyPreperiodic
% <<< [agent-claude 2026-10-29 lean-xref]
Let $\psi$ be coordinatewise squaring on $\PP^2$ and put
\[
L([x:y:z])=[x-y+z:x+y-z:-x+y+z],
\qquad \varphi=\psi\circ L.
\]
This map is attributed to Dupont in
\cite[Example~9]{BIJJLMRS17}. Since $\det L=4$, its critical lines
are $L^{-1}(H_i)$ and map to the coordinate lines $H_i$.
Put $M_0=V(x-y)$, $M_1=V(y-z)$, and $M_2=V(x-z)$.
For example, $\varphi(H_0)=M_0$, since on $x=0$ its first two
coordinates are equal. The other two coordinate lines map to
$M_1$ and $M_2$; direct substitution gives
\[
\varphi(M_0)=M_2,\qquad
\varphi(M_2)=M_1,\qquad
\varphi(M_1)=M_0.
\]
Each restricted map is nonconstant, so these are equalities of
integral effective divisors. Hence the postcritical set consists
of the three $H_i$ and three $M_i$. The branch lines are strictly
preperiodic, while the $M_i$ form a three-cycle.
\end{example}

The normalization in Definition~\ref{def:crit-branch} is essential.
Replacing the normalized branch locus by either the nonreduced-fibre
locus or the non-\'etale locus invalidates the descent theorem.

\begin{example}[Failure for unnormalized branch loci]
\label{ex:pinched}
% >>> [agent-claude 2026-10-29 lean-xref]
% Lean: Pcf.Examples.CharZeroPaperPinched.ambient_three_branch_loci,
%   Pcf.Examples.CharZeroPaperPinched.line_isNormalizedPCF,
%   Pcf.Examples.CharZeroPaperPinched.fiberQ_iso,
%   Pcf.Examples.CharZeroPaperPinched.not_pcf_nonreducedFibre,
%   Pcf.Examples.CharZeroPaperPinched.not_pcf_nonEtale,
%   Pcf.Examples.CharZeroPaperPinched.gScheme_eq_invariantClosureMap
% <<< [agent-claude 2026-10-29 lean-xref]
Let $k$ have characteristic zero and set
\[
R=\{h\in k[z,w]:h(2,1)=h(-2,1)\},\qquad X=\Spec R.
\]
The map $(z,w)\mapsto(z^2,w^2)$ descends to a finite self-map
$f$ of $X$. The normalization $\nu\colon\AA^2\to X$ is an
isomorphism away from the glued point. At the generic point of
every divisor, $f$ therefore agrees with the squaring map, whose
Jacobian is $4zw$. Thus its normalized, nonreduced-fibre, and
non-\'etale \emph{divisorial} branch loci all consist of the
images of $V(z)$ and $V(w)$, both invariant under $f$.
In particular, $f$ is PCF for all three conventions; isolated
non-\'etale points are not components of these branch loci.
The image $B$ of $V(w-1)$ is invariant, and the normalized lift
of $g=\restr{f}{B}$ is $z\mapsto z^2$, so $g$ is PCF for the
normalized branch locus.
However, if $Q\in B$ is represented by $(4,1)$, then the fibre
$g^{-1}(Q)$ is nonreduced. To see this, put $b=z^2$ and
$d=z(z^2-4)$; then
$k[B]=k[b,d]/(d^2-b(b-4)^2)$, and $g^*$ sends
$b\mapsto b^2$ and $d\mapsto b^3-4b$.
The inverse image of $Q=(b,d)=(16,48)$ is the pinched point
$(b,d)=(4,0)$, with fibre $k[d]/(d^2)$.
The forward orbit of $Q$ is the sequence of distinct points
represented by $(4^{2^n},1)$, $n\geq0$.
% >>> [agent-codex 2026-11-04 public-references]
Thus $g$ is not PCF for either unnormalized locus.
% <<< [agent-codex 2026-11-04 public-references]
\end{example}
% <<< [agent-codex 2026-10-27 char0-examples]

%% -------------------------------------------------------
\section{Preliminary lemmas}
\label{sec:prelim}
%% -------------------------------------------------------

\begin{proposition}
	\label{prop:pcf-equiv}
	% >>> [agent-claude 2026-10-29 lean-xref]
	% Lean: Pcf.PaperCharZero.scheme_isNormalizedPCF_iff_branch_preperiodic,
	%   Pcf.PaperCharZero.scheme_isNormalizedPCF_iff_critical_support_preperiodic
	% <<< [agent-claude 2026-10-29 lean-xref]
	Let $f \colon X \to X$ be a finite surjective morphism of a variety. The
	following are equivalent.
	\begin{enumerate}
		\item \label{item:pcf} $f$ is PCF.
		\item \label{item:orbits} Every $D \in \Irr(\branch(f))$ is preperiodic.
		\item \label{item:crit-preperiodic} The critical locus is preperiodic:
		there are $r \geq 0$ and $s \geq 1$ with
		$f^{r+s}(\crit(f)) = f^{r}(\crit(f))$.
	\end{enumerate}
\end{proposition}

\begin{proof}
	Write $\Irr(\branch(f)) = \{D_{1}, \ldots, D_{t}\}$. If $t = 0$ all three
	conditions hold, with $r = 1$ and $s = 1$ in~\eqref{item:crit-preperiodic},
	since then $f(\crit(f)) = \emptyset$; so assume $t \geq 1$.

	\emph{\eqref{item:pcf} implies \eqref{item:orbits}.} The orbit
	$\{f^{n}(D_{i})\}_{n \geq 0}$ lies in the finite set $\PC(f)$, so two of
	its terms coincide, say $f^{r_{i}+s_{i}}(D_{i}) = f^{r_{i}}(D_{i})$ with
	$s_{i} \geq 1$.

	\emph{\eqref{item:orbits} implies \eqref{item:crit-preperiodic}.} Choose
	$r_{i}, s_{i}$ as above and put $r = 1 + \max_{i} r_{i}$ and
	$s = \operatorname{lcm}_{i} s_{i}$. For $n \geq r_{i}$ the relation
	$f^{n+s_{i}}(D_{i}) = f^{n}(D_{i})$ holds, hence so does
	$f^{n+s}(D_{i}) = f^{n}(D_{i})$. Since $f(\crit(f)) = D_{1} \cup \cdots
	\cup D_{t}$,
	\[
	f^{r+s}(\crit(f)) = \bigcup_{i} f^{r-1+s}(D_{i})
	= \bigcup_{i} f^{r-1}(D_{i}) = f^{r}(\crit(f)) .
	\]

	\emph{\eqref{item:crit-preperiodic} implies \eqref{item:pcf}.} Applying
	$f^{j}$ to the relation gives $f^{n+s}(\crit(f)) = f^{n}(\crit(f))$ for
	all $n \geq r$, so every $f^{n}(\crit(f))$ with $n \geq 1$ equals one of
	% >>> [agent-codex 2026-10-28 pcf-periodic-endpoint]
	% $f^{1}(\crit(f)), \ldots, f^{r+s-1}(\crit(f))$.
	$f^{1}(\crit(f)), \ldots, f^{R+s-1}(\crit(f))$, where
	$R=\max\{r,1\}$. This includes $r=0$, $s=1$.
	% <<< [agent-codex 2026-10-28 pcf-periodic-endpoint]
	These are finitely many
	finite unions of integral effective divisors, and $\PC(f)$ is the set of their
	components.
\end{proof}

\begin{lemma}[Branch loci in towers]
	\label{lem:tower}
	% >>> [agent-claude 2026-10-29 lean-xref]
	% Lean: Pcf.PaperCharZero.normalized_branch_tower
	% <<< [agent-claude 2026-10-29 lean-xref]
	Let $f \colon X \to Y$ and $g \colon Y \to Z$ be finite surjective
	morphisms of varieties. Then
	\[
	\branch(g \circ f) = \branch(g) \cup g(\branch(f)) .
	\]
\end{lemma}

\begin{proof}
	By uniqueness of lifts \stacks[(4)]{035Q} the lift of $g \circ f$ to the
	normalizations is $\widetilde{g} \circ \widetilde{f}$. For an integral
	effective divisor
	$\widetilde{P}$ of $\widetilde{X}$ put
	$\widetilde{Q} = \widetilde{f}(\widetilde{P})$ and
	$\widetilde{R} = \widetilde{g}(\widetilde{Q})$. Ramification indices of
	extensions of discrete valuation rings multiply in towers, so
	\[
	e(\widetilde{P}/\widetilde{R})
	= e(\widetilde{P}/\widetilde{Q})\, e(\widetilde{Q}/\widetilde{R}) .
	\]
	Thus $\widetilde{P}$ is critical for $g \circ f$ exactly when it is
	critical for $f$ or $\widetilde{Q}$ is critical for $g$. In the first
	case $\nu_{Z}(\widetilde{R}) = g(\nu_{Y}(\widetilde{Q}))$ lies in
	$g(\branch(f))$, and in the second it lies in $\branch(g)$. This proves
	the inclusion $\subseteq$. Conversely, if $\widetilde{Q}$ is critical for
	$g$, Lemma~\ref{lem:lift} gives an integral effective divisor $\widetilde{P}$
	over it, and the
	displayed formula shows that $\nu_{Z}(\widetilde{g}(\widetilde{Q}))$ lies
	in $\branch(g \circ f)$; and every component of $g(\branch(f))$ has the
	form $\nu_{Z}(\widetilde{R})$ for some $\widetilde{P}$ critical for $f$,
	by~\eqref{eq:branch-normalization}.
\end{proof}

The lemma reflects the picture of a tower $X \to Y \to Z$: a divisor of
$Z$ is branched for the composite when it is branched at the bottom, or
when it is the image of a divisor that is branched one level up.

\begin{proposition}[Iterates]
	\label{prop:iterates}
	% >>> [agent-claude 2026-10-29 lean-xref]
	% Lean: Pcf.PaperCharZero.normalized_iterates
	% <<< [agent-claude 2026-10-29 lean-xref]
	Let $f \colon X \to X$ be a finite surjective morphism of a variety and
	let $m \geq 1$. Then:
	\begin{enumerate}
		\item \label{item:branch-iterate}
		$\branch(f^{m}) = \bigcup_{i=0}^{m-1} f^{i}(\branch(f))$;
		\item \label{item:orbit-propagation} if
		$D \in \Irr(\branch(f))$ then
		$f^{n}(D) \in \Irr(\branch(f^{n+1}))$ for every $n \geq 0$;
		\item \label{item:pc-iterate} $\PC(f^{m}) = \PC(f)$, and
		$\Irr(\branch(f^{m})) \subseteq \PC(f)$;
		\item \label{item:pcf-iterate} $f$ is PCF if and only if $f^{m}$ is PCF.
	\end{enumerate}
\end{proposition}

\begin{proof}[Proof of \eqref{item:branch-iterate}]
	Induction on $m$, using Lemma~\ref{lem:tower} for
	$f^{m+1} = f \circ f^{m}$:
	$\branch(f^{m+1}) = \branch(f) \cup f(\branch(f^{m}))$.
\end{proof}

\begin{proof}[Proof of \eqref{item:orbit-propagation}]
	By~\eqref{item:branch-iterate} with $m = n+1$, the integral effective divisor
	$f^{n}(D)$ is
	contained in $\branch(f^{n+1})$, which is a finite union of integral effective
	divisors, so it
	is one of its components.
\end{proof}

\begin{proof}[Proof of \eqref{item:pc-iterate}]
	By~\eqref{item:branch-iterate} a component of $\branch(f^{m})$ has the
	form $f^{i}(D)$ with $D \in \Irr(\branch(f))$, so
	$\PC(f^{m}) \subseteq \PC(f)$, and in particular
	$\Irr(\branch(f^{m})) \subseteq \PC(f)$. Conversely, for
	$D \in \Irr(\branch(f))$ and $n = qm + i$ with $0 \leq i < m$, the integral
	effective divisor
	$f^{i}(D)$ is a component of $\branch(f^{m})$, and
	$f^{n}(D) = (f^{m})^{q}(f^{i}(D)) \in \PC(f^{m})$.
\end{proof}

\begin{proof}[Proof of \eqref{item:pcf-iterate}]
	This is immediate from~\eqref{item:pc-iterate}.
\end{proof}

\begin{lemma}[PCF under normalization]
	\label{lem:normalization}
	% >>> [agent-claude 2026-10-29 lean-xref]
	% Lean: Pcf.PaperCharZero.normalization_postcritical_and_pcf
	% <<< [agent-claude 2026-10-29 lean-xref]
	Let $f \colon X \to X$ be a finite surjective morphism of a variety and
	let $\widetilde{f} \colon \widetilde{X} \to \widetilde{X}$ be its lift to
	the normalization. Then $\PC(f) = \{\nu(P) : P \in \PC(\widetilde{f})\}$,
	and $f$ is PCF if and only if $\widetilde{f}$ is PCF.
\end{lemma}

\begin{proof}
	By~\eqref{eq:branch-normalization}, $\branch(f) = \nu(\branch(\widetilde{f}))$,
	and $\nu \circ \widetilde{f}^{n} = f^{n} \circ \nu$ for all $n$, so
	$f^{n}(\branch(f)) = \nu(\widetilde{f}^{n}(\branch(\widetilde{f})))$. This
	gives the description of $\PC(f)$, so $f$ is PCF if $\widetilde{f}$ is.
	Conversely, every member of $\PC(\widetilde{f})$ maps into a member of
	$\PC(f)$, and by Lemma~\ref{lem:lift} only finitely many integral effective
	divisors of
	$\widetilde{X}$ map into a given integral effective divisor of $X$.
\end{proof}

\begin{lemma}[Ramified covers]
	\label{lem:covers}
	% >>> [agent-claude 2026-10-29 lean-xref]
	% Lean: Pcf.isNormalizedPCFOfFiniteType_of_cover
	% <<< [agent-claude 2026-10-29 lean-xref]
	Let $f \colon X \to X$, $g \colon Y \to Y$ and $\pi \colon Y \to X$ be
	finite surjective morphisms of varieties with $f \circ \pi = \pi \circ g$.
	If $g$ is PCF, then $f$ is PCF.
\end{lemma}

\begin{proof}
	Put $T = \pi(\PC(g)) \cup \branch(\pi)$, a finite union of integral effective
	divisors of $X$.
	Applying Lemma~\ref{lem:tower} to both sides of $f \circ \pi = \pi \circ g$
	gives
	\[
	\branch(f) \cup f(\branch(\pi)) = \branch(\pi) \cup \pi(\branch(g))
	\subseteq T .
	\]
	So $\branch(f) \subseteq T$, and $T$ is forward invariant: $f(\branch(\pi))
	\subseteq T$ by the display, and
	$f(\pi(\PC(g))) = \pi(g(\PC(g))) \subseteq \pi(\PC(g))$. Hence every
	member of $\PC(f)$ is an integral effective divisor contained in $T$, that is,
	a component of
	$T$, and there are finitely many of these.
\end{proof}

%% -------------------------------------------------------
\section{The bounding theorem}
\label{sec:bounding}
%% -------------------------------------------------------

We bound the branch locus of the restriction of a finite map to an integral
effective divisor.
The two inputs are classical, and we state them in the form used.

\begin{theorem}[Zariski--Nagata purity, \stacks{0BMB}]
\label{thm:purity}
% >>> [agent-claude 2026-10-29 lean-xref]
% Lean: Pcf.PaperCharZero.zariski_nagata_purity
% <<< [agent-claude 2026-10-29 lean-xref]
Let $h \colon X' \to Y'$ be a morphism of locally Noetherian schemes, let
$x \in X'$ and $y = h(x)$. Assume that $\cO_{X',x}$ is normal, that
$\cO_{Y',y}$ is regular, that $h$ is quasi-finite at $x$, that
$\dim \cO_{X',x} = \dim \cO_{Y',y} \geq 1$, and that $h$ is unramified at
every point $x'$ specializing to $x$ with $\dim \cO_{X',x'} = 1$. Then
$h$ is \'etale at $x$.
\end{theorem}

Let $R$ be a Noetherian ring, $a \in R$ a nonzerodivisor with $R/aR$
reduced, $u \in R^{\times}$, and $e \geq 1$ an integer invertible in $R$.
The morphism $\Spec R[z]/(z^{e} - ua) \to \Spec R$ is a \emph{standard
tamely ramified morphism}; for $u = 1$ this is the morphism of
\stacks{0EYF}, and the unit changes nothing below. It is finite locally
free and
\'etale over $\Spec R_{a}$.

\begin{theorem}[Abhyankar's lemma]
\label{thm:abhyankar}
% >>> [agent-claude 2026-10-29 lean-xref]
% Lean: Pcf.PaperCharZero.abhyankar_lemma
% <<< [agent-claude 2026-10-29 lean-xref]
Let $U$ be a locally Noetherian scheme, $D \subset U$ an effective
Cartier divisor that is a regular scheme, and $W \to U \smallsetminus D$
a finite \'etale morphism that is tamely ramified over $U$ in codimension
one. Let $V$ be the normalization of $U$ in $W$. Then \'etale locally on
$U$ the morphism $V \to U$ is a finite disjoint union of standard tamely
ramified morphisms.
\end{theorem}

This is \stacks{0EYG} together with \stacks{0EYH}.

Explicitly, every point of $U$ is in the image of an \'etale morphism
$U' = \Spec R \to U$ such that $D \times_{U} U' = V(a)$ with $a \in R$ a
nonzerodivisor and $R/aR$ regular, and
\[
V \times_{U} U' \cong \coprod_{i \in I}
\Spec R[z_{i}]/(z_{i}^{e_{i}} - u_{i}a)
\]
for a finite set $I$, integers $e_{i} \geq 1$ invertible in $R$ and units
$u_{i} \in R^{\times}$. Over a field of characteristic zero every finite
\'etale cover is tamely ramified in codimension one, since the residue
fields have characteristic zero.

The mechanism of the proof is that a standard piece becomes an isomorphism
when restricted to the divisor. Indeed, for
$C = R[z]/(z^{e} - ua)$ as above,
\begin{equation}\label{eq:standard-piece}
% >>> [agent-claude 2026-10-29 lean-xref]
% Lean: Pcf.PaperCharZero.standard_piece,
%   Pcf.PaperCharZero.standard_piece_canonical
% <<< [agent-claude 2026-10-29 lean-xref]
C/aC = (R/aR)[z]/(z^{e}), \qquad
(C/aC)_{\mathrm{red}} = C/zC = R/aR ,
\end{equation}
where the second equality holds because the regular ring $R/aR$ is
reduced \stacks{00NP}, so that $z$ generates the nilradical of $C/aC$. The
reduced preimage of $V(a)$ in $\Spec C$ maps isomorphically onto $V(a)$.

\begin{lemma}[Integral effective divisors in a proper closed subset]
\label{lem:finitely-many}
% >>> [agent-claude 2026-10-29 lean-xref]
% Lean: Pcf.PaperCharZero.finite_integral_effective_divisors_in_closed
% <<< [agent-claude 2026-10-29 lean-xref]
Let $B$ be a variety and $\Sigma \subsetneq B$ a proper closed subset.
Then $\Sigma$ contains only finitely many integral effective divisors of $B$.
\end{lemma}

\begin{proof}
Each irreducible component of $\Sigma$ has dimension less than
$\dim B$, and an integral effective divisor of $B$ contained in $\Sigma$ has
dimension
$\dim B - 1$, so it is one of the finitely many components of $\Sigma$ of
that dimension.
\end{proof}

For the rest of this section, let $f \colon X \to Y$ be a finite
surjective morphism of normal varieties over a field of characteristic
zero, let $A$ be an integral effective divisor of $X$, let $B = f(A)$, which is
an integral effective divisor of $Y$
by Lemma~\ref{lem:restriction}, and let $g = \restr{f}{A} \colon A \to B$.
Let $Z$ be the union of the components of $\branch(f)$ other than $B$, and
put
\[
\Sigma = (B \cap Z) \cup (B \cap \Sing(Y)) \cup \Sing(B) .
\]

\begin{theorem}[Bounding theorem]
\label{thm:bounding}
% >>> [agent-claude 2026-10-29 lean-xref]
% Lean: Pcf.PaperCharZero.bounding_theorem
% <<< [agent-claude 2026-10-29 lean-xref]
%% >>> [agent-claude 2026-10-04 branch-loci]
With $f$, $A$, $B$, $g$, $Z$ and $\Sigma$ as above, $\Sigma$ is a proper
closed subset of $B$, and $g$ is \'etale at the generic point of every
integral effective divisor of $A$ whose image is not contained in $\Sigma$.
%% <<< [agent-claude 2026-10-04 branch-loci]
\end{theorem}

\begin{proof}
\emph{Step 0: $\Sigma$ is proper.} The set $Z$ is a finite union of
integral effective divisors different from $B$, so $B \not\subseteq Z$. As $Y$
is normal,
$\Sing(Y)$ has codimension at least two, so $B \not\subseteq \Sing(Y)$.
The regular locus of $B$ is open \stacks{07R5} and contains the generic
point, so $\Sing(B) \neq B$. The set $\Sigma$ is closed and, $B$ being irreducible,
proper.

Let $E$ be an integral effective divisor of $B$ not contained in $\Sigma$, with
generic point
$\eta$, so that $\eta \notin \Sigma$. We show that $g$ is \'etale at the
generic point of every integral effective divisor of $A$ over $E$. As
$\Sing(Y)$, $Z$ and
$\Sing(B)$ are closed, there is
an affine open $U \subseteq Y$ containing $\eta$ that is disjoint from
them; then $U$ is regular, $D := B \cap U$ is regular, and $D$ is an
effective Cartier divisor on $U$, because it is integral of codimension
one and the local rings of the regular scheme $U$ are unique
factorization domains \stacks{0AG0}, \stacks{0AGA}. The branch locus of
$f$ over $U$ is $\branch(f) \cap U$, since ramification indices are
computed in local rings, and
\[
\branch(f) \cap U \subseteq (Z \cup B) \cap U = D .
\]
Put $V = f^{-1}(U)$, a normal scheme, finite and surjective over $U$, and
let $f_{V} \colon V \to U$ be the restriction of $f$.

\emph{Step 1: $f_{V}$ is \'etale over $U \smallsetminus D$.} Let
$w \in V \smallsetminus f^{-1}(D)$. If $w$ is the generic point of $V$,
then $f_{V}$ is \'etale at $w$ because the extension of function fields is
separable. Otherwise $\dim \cO_{V,w} = \dim \cO_{U,f(w)} \geq 1$,
because $f$ is finite and surjective between varieties, and we check the
hypotheses of Theorem~\ref{thm:purity}: $\cO_{V,w}$ is normal,
$\cO_{U,f(w)}$ is regular, and $f_{V}$ is quasi-finite. A codimension-one
point $w'$ specializing to $w$ does not lie in the closed set
$f^{-1}(D)$, so $f(w') \notin \branch(f)$; hence $e = 1$ at $w'$, and
$f_{V}$ is \'etale, in particular unramified, at $w'$ by
Lemma~\ref{lem:dictionary}. Purity shows that $f_{V}$ is \'etale at $w$.

\emph{Step 2: Abhyankar's lemma.} By Step~1,
$W := V \smallsetminus f^{-1}(D) \to U \smallsetminus D$ is finite
\'etale, and it is tamely ramified over $U$ in codimension one because
the characteristic is zero. Moreover $V$ is the normalization of $U$ in
$W$: the normalization of $U$ in $W$ factors through the finite morphism
$V \to U$ \stacks[(2)]{035I}, and the normal ring of an affine open of $V$
is integrally closed in the ring of functions on its dense open
intersection with $W$ \stacks{030C}. Choose an \'etale
$\rho \colon U' = \Spec R \to U$ with $\eta$ in its image as in
Theorem~\ref{thm:abhyankar}, so that
\[
V' := V \times_{U} U' \cong \coprod_{i \in I}
\Spec R[z_{i}]/(z_{i}^{e_{i}} - u_{i}a), \qquad D' := D \times_{U} U' =
V(a) .
\]

\emph{Step 3: restriction to the divisor.} The closed set
$f^{-1}(D) \subseteq V$ has the same dimension as $D$, and $A \cap V$ is
one of its irreducible components. Its preimage $T$ in $V'$ is reduced,
because $V' \to V$ is \'etale, and is a closed subset of the reduced
preimage of $D'$, which by~\eqref{eq:standard-piece} is the disjoint
union of copies $V(z_{i})$ of the regular scheme $D'$. The \'etale, hence
open, map $T \to A \cap V$ sends the generic point of each component of
$T$ to the generic point of the irreducible $A \cap V$, so $T$ is
equidimensional of dimension $\dim D$; and $D'$, being \'etale over the
irreducible $D$, is equidimensional of the same dimension. A reduced closed
subset of full dimension in the regular scheme $\coprod_{i} V(z_{i})$,
whose connected components are integral \stacks{00NP}, is a union of
connected components. So $T$ is a union of connected components of
$\coprod_{i} V(z_{i})$, each mapping isomorphically onto a connected
component of $D'$ under the base change $f' = f \times_{U} U'$.

\emph{Step 4: conclusion.} Let $\xi \in A$ be a point with
$g(\xi) = \eta$. Its closure maps finitely onto $E$, so $\xi$ is the
generic point of an integral effective divisor of $A$ over $E$, and every
integral effective divisor of $A$ over $E$
arises this way. The point $\xi$ lies in $A \cap V$, and we choose a point
$\xi'$ of $T$ over it. By Step~3 the morphism $T \to D'$ is \'etale at $\xi'$,
and
$D' \to D$ is \'etale. The composites $T \to A \cap V \to D$ and
$T \to D' \to D$ agree, both being induced by $f'$, and $T \to A \cap V$
is \'etale with image an
open neighbourhood of $\xi$; since \'etaleness is \'etale local on the
source \stacks{036W}, the morphism $A \cap V \to D$ is \'etale at $\xi$.
%% >>> [agent-claude 2026-10-04 branch-loci]
As $A \cap V$ and $D$ are open in $A$ and $B$, the map $g$ is \'etale at
$\xi$. As $\xi$ was an arbitrary point over $\eta$, this proves the
claim.
%% <<< [agent-claude 2026-10-04 branch-loci]
\end{proof}

By Lemma~\ref{lem:finitely-many}, $\Sigma$ contains only finitely many
integral effective divisors of $B$. This, and not the inclusion alone, is what
the dynamical
argument uses.
No hypothesis relating $B$ to
$\branch(f)$ is imposed, and none would help. If $B$ is not a branch
component, then $\branch(f) \cap U = \emptyset$, purity alone makes $V$
finite \'etale over $U$, and every $e_{i}$ equals one.

% >>> [agent-codex 2026-10-28 branch-bound-examples]
% It is natural to expect that the restriction ramifies exactly where $B$
% meets the other branch components, that is, $\branch(g) = B \cap Z$
% away from the singular loci. The inclusion is all one can have: the
% first example below shows that the equality fails for the normalized
% branch locus, and the second that
% the term $B \cap \Sing(Y)$ cannot be omitted.
The first example disproves the global equality
$\branch(g)=B\cap Z$ at a point of $\Sing(B)$. The second disproves
the equality even where $B$ and $Y$ are smooth. The remaining
examples show why the two singular-locus terms occur in the bound.
% <<< [agent-codex 2026-10-28 branch-bound-examples]

\begin{example}[The equality $\branch(g) = B \cap Z$ fails]
\label{ex:cusp}
% >>> [agent-claude 2026-10-29 lean-xref]
% Lean: Pcf.Examples.CharZeroPaperCuspCanonical.normalized_branch_strict_inclusion,
%   Pcf.Examples.CharZeroPaperCuspCanonical.cusp_inter_normalized_ambient_branch,
%   Pcf.Examples.CharZeroPaperCuspCanonical.cuspOrigin_not_regular,
%   Pcf.Examples.CharZeroPaperCuspSourceLine.sourceLine_isClosedImmersion_range
% <<< [agent-claude 2026-10-29 lean-xref]
Let $f \colon \AA^{2} \to \AA^{2}$ be $f(u, v) = (u^{2}, u^{3} + v)$ and
$A = V(v)$. Then $g(u) = (u^{2}, u^{3})$ maps $A$ onto the cuspidal cubic
$B = V(y^{2} - x^{3})$. The Jacobian determinant of $f$ is $2u$, so
$\branch(f) = V(x)$, which is not $B$; thus $Z = V(x)$ and
$B \cap Z = \{(0, 0)\}$. But $g$ is the normalization of $B$, so its lift
to normalizations is the identity of $\AA^{1}$, and
$\branch(g) = \emptyset \subsetneq B \cap Z$.
\end{example}

Here $(0,0) \in \Sing(B)$, so Theorem~\ref{thm:bounding} makes no claim at
that point.

% >>> [agent-codex 2026-10-28 branch-bound-examples]
\begin{example}[Strict inclusion on the regular locus]
\label{ex:smooth-restriction}
% >>> [agent-claude 2026-10-29 lean-xref]
% Lean: Pcf.Examples.CharZeroPaperSmoothBranch.normalized_branch_strict_inclusion,
%   Pcf.Examples.CharZeroPaperSmoothBranch.parabola_inter_normalized_ambient_branch,
%   Pcf.Examples.CharZeroPaperSmoothBranch.tangent_ideals_equal
% <<< [agent-claude 2026-10-29 lean-xref]
Over $\mathbb C$ let
$f\colon X=\AA^2_{u,v}\to Y=\AA^2_{x,y}$ be
$f(u,v)=(u^2,v)$, and take $A=V(v-u)$.
The coordinate ring $\mathbb C[u,v]$ is free of rank two over
$\mathbb C[u^2,v]$, so $f$ is finite and surjective.
Its Jacobian is $2u$, hence $\branch(f)=V(x)$.
The image $B=f(A)=V(x-y^2)$ is a smooth integral effective divisor,
and $Z=V(x)$ is the sole branch component of $f$.
On $A$, with parameter $t=u=v$, the restriction is
$g(t)=(t^2,t)$, an isomorphism onto $B$ with inverse $y$.
Thus its normalized branch locus is empty, while
\[
\branch(g)=\emptyset\subsetneq B\cap Z=\{(0,0)\}.
\]
Both $B$ and $Y$ are smooth at the failure point. The curves $B$
and $Z$ are tangent there; no transversality hypothesis occurs in
Theorem~\ref{thm:bounding}.
\end{example}
% <<< [agent-codex 2026-10-28 branch-bound-examples]

\begin{example}[An ambient singularity contributes a branch point]
\label{ex:branch-singular}
% >>> [agent-claude 2026-10-29 lean-xref]
% Lean: Pcf.Examples.CharZeroPaperConeBranch.ex_branch_singular
% <<< [agent-claude 2026-10-29 lean-xref]
Let $Y = \Spec k[x, y, z]/(xz - y^{2})$, a normal surface singular only
at its vertex $o$, and let $f \colon \AA^{2} \to Y$ be
$f(u, v) = (u^{2}, uv, v^{2})$, the quotient by $(u, v) \mapsto (-u, -v)$.
It is finite and surjective, and \'etale away from the origin, which is
the only point over $o$; so $\branch(f) = \emptyset$ and $Z = \emptyset$.
For $A = V(v)$ and $B = V(y, z)$ the restriction is $g(u) = u^{2}$ from
one affine line to another, and $\branch(g) = \{o\}$, while
$B \cap Z = \emptyset$.
\end{example}

% >>> [agent-codex 2026-10-28 branch-bound-examples]
% The third term is in a different position. It is used in the proof,
% through the regularity of $D$, but we do not know an example that forces
% it.
%
% \begin{question}
% \label{q:sing-B}
% Is there a finite surjective morphism $f \colon X \to Y$ of normal
% varieties over a field of characteristic zero, with integral effective
% divisors $A$ of $X$ and $B = f(A)$ of $Y$, such that
% $\branch(\restr{f}{A}) \not\subseteq
% (B \cap Z) \cup (B \cap \Sing(Y))$?
% \end{question}
%
% We queried this question with an Astra and the machine was unable
% to come up with an example.

\begin{example}[A singular image divisor contributes a branch point]
\label{ex:branch-sing-B}
% >>> [agent-claude 2026-10-29 lean-xref]
% Lean: Pcf.Examples.CharZeroPaperSingularImageGeometry.sing_term_necessary,
%   Pcf.Examples.CharZeroPaperSingularImageBounding.singB_term_cannot_be_dropped,
%   Pcf.Examples.CharZeroPaperSingularImageNormalizedBranch.planeMap_normalizedBranch
% <<< [agent-claude 2026-10-29 lean-xref]
Work over $\mathbb{C}$, choose $a \in \mathbb{C}$ with $a^{2}=-3$,
and put
\[
F=x^{4}+2ax^{2}y^{2}+y^{4},\qquad
G=x^{5}y-xy^{5},\qquad
H=x^{4}-2ax^{2}y^{2}+y^{4}.
\]
Let $f=(F,G)\colon X=\AA^{2}\to Y=\AA^{2}$, with target coordinates
$U,V$, and let $B=V(U^{3}-12aV^{2})$. Direct expansion gives
\[
\det Df=4H^{2},\qquad F^{3}-12aG^{2}=H^{3}.
\]
The second identity makes $H$ integral over $\mathbb{C}[F,G]$.
Both $x^{4}$ and $y^{4}$ satisfy
\[
T^{2}-\frac{F+H}{2}T-\frac{(F-H)^{2}}{48}=0,
\]
so transitivity of integrality makes $x,y$ integral over
$\mathbb{C}[F,G]$, and $f$ is finite. Its nonzero Jacobian gives
dominance; hence $f$ is surjective.

The polynomial $h(z)=z^{4}-2az^{2}+1$ has four distinct roots:
$h(0)=1$, while $h(z)=4$ when $z^{2}=a$, so it has no common root
with $h'(z)=4z(z^{2}-a)$. Thus $H=0$ is a union of four distinct
lines $A_{\lambda}=V(y-\lambda x)$. For each root $\lambda$ of $h$,
the restriction to $A_{\lambda}$ is
\[
t\longmapsto(c_{4}t^{4},c_{6}t^{6}),\qquad
c_{4}=4a\lambda^{2},\quad c_{6}=\lambda(1-\lambda^{4}),
\]
with $c_{4}c_{6}\ne0$ and $c_{4}^{3}=12ac_{6}^{2}$.
Indeed $\lambda\ne0$, and $\lambda^{4}=1$ would imply $a^{2}=1$
from $h(\lambda)=0$. The irreducible cuspidal curve $B$ is the
image of every $A_{\lambda}$. Each image is a nonconstant closed
curve in $B$, and $U^{3}/(12a)$ is not a square in
$\mathbb{C}(U)$, so $B$ is integral.
Since $f^{*}(U^{3}-12aV^{2})=H^{3}$, the ramification index of
each $A_{\lambda}$ over $B$ is three. The Jacobian is invertible
away from $H=0$, so there are no other critical integral
effective divisors. Consequently $\branch(f)=B$ and $Z=\emptyset$.

Take $A=A_{\lambda}$. The normalization of $B$ is
$s\mapsto(c_{4}s^{2},c_{6}s^{3})$: its coordinate ring is
$\mathbb{C}[s^{2},s^{3}]\subset\mathbb{C}[s]$, and
$s=(c_{4}/c_{6})V/U$ lies in its function field. The normalized
lift of $g=\restr{f}{A}$ is $t\mapsto t^{2}$. Hence
\[
\branch(g)=\{(0,0)\},\qquad
(B\cap Z)\cup(B\cap\Sing(Y))=\emptyset.
\]
The branch point lies in $\Sing(B)$, showing that the third term
of Theorem~\ref{thm:bounding} cannot in general be omitted.
\end{example}

\begin{remark}[Sharpness questions for the branch-locus bound]
\label{rem:branch-bound-sharpness}
% >>> [agent-claude 2026-10-29 lean-xref]
% Lean: Pcf.Examples.CharZeroPaperSingularImageBounding.singB_term_cannot_be_dropped,
%   Pcf.Examples.CharZeroPaperConeBranch.ex_branch_singular,
%   Pcf.Examples.CharZeroPaperSmoothBranch.normalized_branch_strict_inclusion
% <<< [agent-claude 2026-10-29 lean-xref]
An earlier draft asked whether $\Sing(B)$ was needed in the
exceptional set of Theorem~\ref{thm:bounding}.
Example~\ref{ex:branch-sing-B} answers this question affirmatively.
Example~\ref{ex:cusp} shows failure of global equality
$\branch(g)=B\cap Z$, while
Example~\ref{ex:smooth-restriction} shows failure even where $B$
and $Y$ are smooth. Example~\ref{ex:branch-singular} shows that
$B\cap\Sing(Y)$ can also contribute a branch point.
\end{remark}
% <<< [agent-codex 2026-10-28 branch-bound-examples]

%% -------------------------------------------------------
\section{Invariant integral effective divisors of normal varieties}
\label{sec:main}
%% -------------------------------------------------------

\begin{theorem}
\label{thm:normal-invariant}
% >>> [agent-claude 2026-10-29 lean-xref]
% Lean: PaperCharZero.theoremA
% <<< [agent-claude 2026-10-29 lean-xref]
Let $X$ be a normal variety over a field of characteristic zero, let
$f \colon X \to X$ be a finite surjective PCF morphism, and let $A$ be an
% >>> [agent-codex 2026-10-28 statement-notation-charzero]
integral effective divisor of $X$ with $f(A) = A$. Then
$\restr{f}{A} \colon A \to A$ is PCF.
% <<< [agent-codex 2026-10-28 statement-notation-charzero]
%% >>> [agent-claude 2026-10-04 branch-loci]
%PCF. %taylor did this
%% <<< [agent-claude 2026-10-04 branch-loci]
\end{theorem}

\begin{proof}
% >>> [agent-codex 2026-10-28 statement-notation-charzero]
Put $g=\restr{f}{A} \colon A\to A$.
% <<< [agent-codex 2026-10-28 statement-notation-charzero]
Let $W_{A}$ be the union of the members of $\PC(f)$ other than $A$, a
closed subset of $X$ because $\PC(f)$ is finite, and put
\[
\Sigma = (A \cap W_{A}) \cup (A \cap \Sing(X)) \cup \Sing(A) .
\]
This is a proper closed subset of $A$, by Step~0 of the proof of
Theorem~\ref{thm:bounding}, so by Lemma~\ref{lem:finitely-many} it
contains only finitely many integral effective divisors of $A$. The point of the
proof is that
$\Sigma$ does not depend on any iterate of $f$.

Fix $m \geq 1$. Then $g^{m} = \restr{f^{m}}{A}$. Let $Z_{m}$ be the union
of the components of $\branch(f^{m})$ other than $A$; by
Proposition~\ref{prop:iterates}\eqref{item:pc-iterate} these components
lie in $\PC(f)$, so $Z_{m} \subseteq W_{A}$. Theorem~\ref{thm:bounding},
applied to $f^{m} \colon X \to X$ and the integral effective divisor $A$, whose
image is $A$,
%% >>> [agent-claude 2026-10-13 consistency]
shows that $g^{m}$ is \'etale at the generic point of every integral
effective divisor of $A$ whose image is not contained in
$(A \cap Z_{m}) \cup (A \cap \Sing(X)) \cup \Sing(A)$.
% >>> [agent-codex 2026-10-27 char0-local-index]
The image generic point lies outside $\Sing(A)$, so
Lemma~\ref{lem:etale-normalized-index} shows that the divisor is not
critical for the normalized branch locus. Thus
% <<< [agent-codex 2026-10-27 char0-local-index]
\[
\branch(g^{m}) \subseteq (A \cap Z_{m}) \cup (A \cap \Sing(X)) \cup
\Sing(A) \subseteq \Sigma .
\]

Now let $D \in \Irr(\branch(g))$. By
Proposition~\ref{prop:iterates}\eqref{item:orbit-propagation}, for each
$n \geq 0$ the integral effective divisor $g^{n}(D)$ is a component of
$\branch(g^{n+1})$, so it is an integral effective divisor of $A$
contained in $\Sigma$. There are finitely many of those, so the orbit of
$D$ is finite. As $\branch(g)$ has finitely many components, $\PC(g)$
is finite.
%% <<< [agent-claude 2026-10-13 consistency]
\end{proof}

%% [agent-claude 2026-10-13 consistency] removed Remark
%% rem:normal-invariant-pcf and the \iffalse'd Corollary
%% cor:normal-invariant-pcf; the 1-PCF form is
%% Proposition~\ref{prop:normal-unnormalized} in Appendix A.
%% <<< [agent-claude 2026-10-04 branch-loci]

\begin{corollary}[Preperiodic integral effective divisors of normal varieties]
\label{cor:normal-preperiodic}
% >>> [agent-claude 2026-10-29 lean-xref]
% Lean: PaperCharZero.theoremA
% <<< [agent-claude 2026-10-29 lean-xref]
Let $X$ be a normal variety over a field of characteristic zero, let
$f \colon X \to X$ be a finite surjective PCF morphism, and let $A$ be an
integral effective divisor of $X$, and let $r \geq 0$ and $s \geq 1$ be integers
with
$f^{r+s}(A) = f^{r}(A)$. Then $\restr{f^{s}}{f^{r}(A)}$ is PCF.
%% >>> [agent-claude 2026-10-04 branch-loci]
%$1$-PCF, hence $0$-PCF and PCF.
%% <<< [agent-claude 2026-10-04 branch-loci]
\end{corollary}

\begin{proof}
The integral effective divisor $B = f^{r}(A)$ satisfies $f^{s}(B) = B$, and $F =
f^{s}$ is PCF by
Proposition~\ref{prop:iterates}\eqref{item:pcf-iterate}. Apply
%% [agent-claude 2026-10-13 consistency] citation of the removed
%% cor:normal-invariant-pcf dropped.
Theorem~\ref{thm:normal-invariant} to $F$ and $B$.
\end{proof}

%% [agent-claude 2026-10-13 consistency] removed, as approved, the remark
%% "Remark rem:normal-invariant-pcf actually give the stronger conclusion
%% that the restriction is $PCF_0$"; see
%% Proposition~\ref{prop:normal-unnormalized}.

%% -------------------------------------------------------
\section{The general case}
\label{sec:general}
%% -------------------------------------------------------

We remove the normality of $X$ by passing to the normalization, and
transport the result back with Lemma~\ref{lem:covers}.

\begin{theorem}
\label{thm:general}
% >>> [agent-claude 2026-10-29 lean-xref]
% Lean: PaperCharZero.theoremA
% <<< [agent-claude 2026-10-29 lean-xref]
Let $X$ be a variety over a field of characteristic zero, let
$f \colon X \to X$ be a finite surjective PCF morphism, and let $A$ be an
integral effective divisor of $X$, and let $r \geq 0$ and $s \geq 1$ be integers
with
$f^{r+s}(A) = f^{r}(A)$. Then $\restr{f^{s}}{f^{r}(A)}$ is PCF.
\end{theorem}

\begin{proof}
Put $B = f^{r}(A)$, an integral effective divisor with $f^{s}(B) = B$, and
$g = \restr{f^{s}}{B}$. Let $\nu \colon \widetilde{X} \to X$ be the
normalization and
$\widetilde{f}$ the lift of $f$, so that $\nu \circ \widetilde{f} = f
\circ \nu$; by Lemma~\ref{lem:normalization}, $\widetilde{f}$ is PCF. Let
$\mathcal{S}$ be the set of integral effective divisors of $\widetilde{X}$ that
map onto $B$; it
is nonempty and finite by Lemma~\ref{lem:lift}. For
$\widetilde{B} \in \mathcal{S}$ the image $\widetilde{f}^{s}(\widetilde{B})$
is an integral effective divisor with
$\nu(\widetilde{f}^{s}(\widetilde{B})) = f^{s}(B) = B$, so $\widetilde{f}^{s}$
maps $\mathcal{S}$ to itself. A self-map of a finite set has a periodic
point, so there are $\widetilde{B} \in \mathcal{S}$ and $c \geq 1$ with
$\widetilde{f}^{sc}(\widetilde{B}) = \widetilde{B}$.

The map $\widetilde{f}^{sc}$ is PCF by
Proposition~\ref{prop:iterates}\eqref{item:pcf-iterate}, and
%% [agent-claude 2026-10-13 consistency] citation updated
$\widetilde{X}$ is normal, so Corollary~\ref{cor:normal-preperiodic}
shows that
$h = \restr{\widetilde{f}^{sc}}{\widetilde{B}}$ is PCF. The restriction
$\pi = \restr{\nu}{\widetilde{B}} \colon \widetilde{B} \to B$ is finite and
surjective by Lemma~\ref{lem:restriction}, and
$\pi \circ h = g^{c} \circ \pi$. By Lemma~\ref{lem:covers}, $g^{c}$ is PCF,
and so is $g$ by Proposition~\ref{prop:iterates}\eqref{item:pcf-iterate}.
\end{proof}

\begin{corollary}
\label{cor:pcfawd}
% >>> [agent-claude 2026-10-29 lean-xref]
% Lean: PaperCharZero.cor_pcfawd
% <<< [agent-claude 2026-10-29 lean-xref]
Conjecture~\ref{conj:pcfawd} holds: every finite surjective PCF self-map
of a variety over a field of characteristic zero is PCF all the way down.
\end{corollary}

\begin{proof}
Induction on $\dim X$, the case $\dim X = 0$ being part of
Definition~\ref{def:pcfawd}. If $\dim X \geq 1$, $f$ is PCF, and $A$, $r$,
$s$, $B$ are as in Definition~\ref{def:pcfawd}, then
$\restr{f^{s}}{B}$ is a finite surjective self-map of the variety $B$
(Lemma~\ref{lem:restriction}), of dimension $\dim X - 1$, and it is PCF by
Theorem~\ref{thm:general}. By induction it is PCF all the way down.
\end{proof}

The descent continues along every chain of periodic integral effective divisors.
It does not
assert that a periodic integral effective divisor exists at each stage; if none
does, the
chain stops there.

\appendix

%% -------------------------------------------------------
% >>> [agent-codex 2026-10-27 char0-appendix-split]
% The detailed alternate-branch-locus comparison now belongs to the
% companion paper. Example~\ref{ex:pinched} retains its negative result
% in this manuscript.
%% \taylor{I copy-pasted this here it may be out of order.} -- addressed
%% \taylor{Why do we have a proposition AND a remark? This is dragging on.}
%% Addressed: the detailed comparison moved to the companion article.
%% \taylor{The above remark needs to be folded into here.}
%% Addressed in the companion appendix.
%% \taylor{Move to appendix} (twice) -- addressed in the companion article.
%% \taylor{for both what?} -- addressed: both unnormalized loci.
% <<< [agent-codex 2026-10-27 char0-appendix-split]

%% -------------------------------------------------------
\section{Software and AI methods}
\label{app:ai}
%% -------------------------------------------------------

% >>> [agent-codex 2026-10-27 char0-ai-detail]

\subsection{Typesetting}
AI agents helped transcribe and revise the LaTeX manuscript using
the mathcity LaTeX workflow
\cite[\href{https://github.com/tdupu/mathcity/tree/0205e0adfb8963a843c28616fcedceafc7ca5680}{public workflow source}]{Mathcity26}. This revision was checked
for line length, stale labels, and citations, and built with
\texttt{latexmk}~\cite[\href{https://ctan.org/pkg/latexmk}{CTAN package page}]{Latexmk26}.
% >>> [agent-codex 2026-11-04 public-references]
A successful build is not a proof check.
% <<< [agent-codex 2026-11-04 public-references]
% >>> [agent-codex 2026-10-28 public-workflow-source]
A pinned public collection of related skills is
\cite[\href{https://github.com/tdupu/agent-skills-public/tree/1e1b9ba1e6bcf20c31f0bacc799ecdc68b502257}{public skill collection}]{AgentSkillsPublic26}.
% <<< [agent-codex 2026-10-28 public-workflow-source]

\subsection{Agent management}
The authors directed the mathematical scope and editorial choices.
Agents worked on bounded drafting, source checks, and review tasks.
The public LaTeX, style, and AI policies
\cite[\href{https://github.com/tdupu/mathcity/tree/0205e0adfb8963a843c28616fcedceafc7ca5680}{public repository policies}]{Mathcity26}
governed the manuscript's preparation. Agent verdicts do
not replace author review.

\subsection{Natural language proofs}
The initial argument of the second and third authors reduced descent
to an equality of branch loci that does not hold. The replacement
upper bound in Theorem~\ref{thm:bounding} was developed from a
frontier-dump suggestion and then expanded into the local proof
given here. Theorem~\ref{thm:general} uses that bound for every
iterate; its proof is part of this manuscript.

\subsection{Adversarial reviews}
AI agents challenged the local purity and Abhyankar steps, the
examples, and the translation into the manuscript.
% >>> [agent-codex 2026-10-27 referee-skill-name]
%% \taylor{referee-report} -- addressed as a typeset skill name below.
The \skill{referee-report} review of the transcription identified
provenance and literature omissions that were repaired in later drafts.
% <<< [agent-codex 2026-10-27 referee-skill-name]

\subsection{AI assisted computations}
The squaring and Dupont calculations in earlier drafts were checked
% >>> [agent-codex 2026-10-28 software-pinpoints-charzero]
in SymPy \cite[\href{https://doi.org/10.7717/peerj-cs.103}{software paper}]{SymPy17},
\href{https://github.com/sympy/sympy/releases/tag/1.13.2}%
{version~1.13.2}.
% >>> [agent-codex 2026-11-04 public-references]
The ring and fibre calculations for
Example~\ref{ex:pinched} were checked in Singular~4.4.1
\cite[\href{https://www.singular.uni-kl.de/}{software site}; version~4.4.1]{Singular}.
% <<< [agent-codex 2026-11-04 public-references]
% <<< [agent-codex 2026-10-28 software-pinpoints-charzero]
Example~\ref{ex:koch} is quoted
from the published Belk--Koch paper.

\subsection{Generation of the AI statement}
An earlier disclosure was assembled with the
\skill{write-materials-and-methods} workflow and checked with
\skill{latex-ai-statement}
\cite[\href{https://github.com/tdupu/mathcity/tree/0205e0adfb8963a843c28616fcedceafc7ca5680}{public workflow source}]{Mathcity26}. OpenAI Codex
\cite[\href{https://openai.com/codex/}{harness page}]{OpenAICodex26}
(GPT-6 family
\cite[\href{https://developers.openai.com/api/docs/guides/latest-model}{model-family guide}]{OpenAIGPT6Family26};
exact deployment identifier and effort not exposed in this session)
drafted this paper-specific disclosure and the short
introduction notice using the \skill{using-latexpowers} and
\skill{math-workflow} instructions. The present text has received
agent review for mathematical dependencies and a LaTeX build;
human review and acceptance of this new wording remain pending.
The authors retain responsibility for its accuracy.
% >>> [agent-codex 2026-10-27 editorial-ai-statement]
% <<< [agent-codex 2026-10-27 editorial-ai-statement]

\subsection{Proofreading and adversarial review}
The earlier transcription received a line-by-line agent referee
report, and the alternate-locus example received a separate
adversarial calculation. The present shortened text was checked
for proof dependencies, references, and LaTeX compilation.
A compiled PDF cannot verify a mathematical assertion.

\subsection{Formalization}\label{sec:formalization}
Theorem~\ref{thm:A} (equivalently Theorem~\ref{thm:general})
and Corollary~\ref{cor:pcfawd} are formalized in Lean~4
% >>> [agent-codex 2026-10-28 mathlib-pin-charzero]
\cite[\href{https://link.springer.com/chapter/10.1007/978-3-030-79876-5_37}{prover paper}]{Lean4},
\href{https://github.com/leanprover/lean4/releases/tag/v4.35.0}%
{version~\texttt{v4.35.0}}, against Mathlib
\cite[\href{https://github.com/leanprover-community/mathlib4/tree/56a005855ea17e0b0958f73ec801ae140b0a1679}{source revision~\texttt{56a005855ea1}}]{MathlibRev26} at tag
% <<< [agent-codex 2026-10-28 mathlib-pin-charzero]
\texttt{v4.35.0}. 

% >>> [agent-codex 2026-11-04 public-references]
The paper-facing declarations are
\texttt{PaperCharZero.theoremA} and
\texttt{PaperCharZero.cor\_pcfawd}.
The first applies the actual return map to the reduced image
of a preperiodic integral effective divisor; the second proves
the divisor recursion at every finite depth. A proved
dimension-bound equivalence includes the automatic dimension-zero
case. 

The Lean setting
allows integral, locally finite type schemes with Noetherian
hypotheses and does not require separatedness, so it includes
the varieties considered here. Its PCF predicate is proved
equivalent to the normalized branch-locus definition used in
this paper. 

The project imports the proved purity and
Abhyankar inputs from its shared formalization library;
they are not additional hypotheses of these two declarations.
The dedicated project built successfully, and both declarations
print only \texttt{propext}, \texttt{Classical.choice}, and
\texttt{Quot.sound} as transitive axioms. The formalization
covers these headline claims.
% >>> [agent-claude 2026-10-28 char0-support-coverage]

Separate imported modules prove Theorem~\ref{thm:bounding} as
stated, with the exceptional set $\Sigma$ of the theorem, together
with Lemmas~\ref{lem:lift}, \ref{lem:dictionary},
\ref{lem:etale-normalized-index}, \ref{lem:restriction},
\ref{lem:tower}, \ref{lem:normalization}, \ref{lem:covers} and
\ref{lem:finitely-many}, Propositions~\ref{prop:pcf-equiv}
and~\ref{prop:iterates}, Theorems~\ref{thm:purity}
and~\ref{thm:abhyankar}, and the
isomorphisms~\eqref{eq:standard-piece}. They also prove every
example: Examples~\ref{ex:squaring-PN}, \ref{ex:koch},
\ref{ex:dupont}, \ref{ex:pinched}, \ref{ex:cusp},
\ref{ex:smooth-restriction}, \ref{ex:branch-singular}
and~\ref{ex:branch-sing-B}. For Example~\ref{ex:branch-sing-B} this
includes the statement that the term $\Sing(B)$ cannot be removed
from $\Sigma$ in Theorem~\ref{thm:bounding}. 

There are some minor differences between the lean and the text which we now enumerate:
In Example~\ref{ex:branch-singular} the quotient by $\pm 1$ is
formalized as the ring of invariants. In Example~\ref{ex:pinched}
the nonreduced-fibre and non-\'etale branch loci are defined for
the example itself. The comparison with the classical condition on
$\PP^{1}$ in the introduction is proved for maps given by
basepoint-free binary forms. For Remark~\ref{rem:why-normalized},
the statement that the normalization of the cusp is an unbranched
homeomorphism is formalized; the complex-analytic local degree is
not.
% <<< [agent-claude 2026-10-28 char0-support-coverage]
The formalization of Theorem~\ref{thm:A} and
Corollary~\ref{cor:pcfawd} is registered at
Palomar~\cite[entry~PALOMAR-2026-10-30-000004 v1]{palomar-2026-10-30-000004-v1}
(\href{https://palomar-registry.org/entry?id=PALOMAR-2026-10-30-000004&version=1}{entry~PALOMAR-2026-10-30-000004 v1}).
The public source repository provides the Lean code and a
statement-to-declaration
map~\cite[characteristic-zero source map]{PcfawdFormalization26}.
% >>> [agent-codex 2026-11-04 public-blueprint-link]
A public \href{https://github.com/tdupu/pcfawd-formalization/tree/main/char-zero/blueprint}{source blueprint}
links the mathematical claims of this article to the corresponding
Lean declarations and records their compiler checks.
% <<< [agent-codex 2026-11-04 public-blueprint-link]
% <<< [agent-codex 2026-11-04 public-references]
% <<< [agent-codex 2026-10-27 char0-ai-detail]

%% -------------------------------------------------------
\subsection{Record keeping and token accounting}
\label{app:tokens}
%% -------------------------------------------------------
% >>> [agent-codex 2026-11-04 public-references]
The authors kept internal AI-usage and token records for the work
summarized in Table~\ref{tab:tokens}.
% <<< [agent-codex 2026-11-04 public-references]

Table~\ref{tab:tokens} lists the AI work whose output is used in this
paper. Counts are provider per-request usage records where available;
otherwise the table labels harness-reported aggregates or final goal
telemetry. Per-request counts are measured at a cutoff, so later requests
in the same session are not included. Tokens are split into uncached input,
cache writes, cache reads and output where the harness reports the split.
Costs are API-equivalent token charges
only, computed as the sum over token classes of tokens times the
provider's published rate, and are not actual billing; per-call charges
such as web searches are excluded. Rates per million tokens: for Claude
Opus~5.5, input $4.00$, cache write $8.00$ (one hour) or $5.00$ (five
minutes), cache read $0.20$, output $20.00$; for GPT-6 Astra, input $10.00$, cache read $1.00$, output $50.00$
at short-context rates, and twice the input and cache-read rates and
$75.00$ for output at long-context rates, which gives the range shown. No
rate for Claude Fable~5.1
\cite[\href{https://www.anthropic.com/claude-fable-and-mythos-5-1}{model page}]{AntFable51}
was available, and the Codex goal telemetry has no
token-class split or applicable rate, so those rows are unpriced, which is
not a claim of zero cost. Earlier sessions whose usage was not recorded are
listed as such.
The present paper-one revision and its bounded agent review have no
exposed per-task token count or reliable aggregate and are unpriced.
% >>> [agent-codex 2026-10-27 editorial-token-note]
The later cross-paper editorial revision also has no exposed
per-task count and is unpriced.
% <<< [agent-codex 2026-10-27 editorial-token-note]
% >>> [agent-codex 2026-10-28 coauthor-feedback-accounting]
The coauthor-feedback row below uses measured per-response counts
through a fixed cutoff; later requests are
excluded. Its charge is a Standard API equivalent, not observed billing.
The GPT-6 Sol rates per million tokens are $2.00$ for uncached input,
$0.20$ for cached input, $2.50$ for cache writes, and $10.00$ for
output
\cite[\href{https://developers.openai.com/api/docs/models/gpt-6-sol}{GPT-6 Sol model page}]{OpenAISol26}
\cite[\href{https://openai.com/codex/}{Codex harness page}]{OpenAICodex26}.
% <<< [agent-codex 2026-10-28 coauthor-feedback-accounting]

\begin{table}[ht]
\centering
\small
\renewcommand{\arraystretch}{1.2}
\begin{tabular}{@{}p{0.27\textwidth}p{0.25\textwidth}%
p{0.27\textwidth}p{0.11\textwidth}@{}}
\hline
Task & Model and harness & Tokens & Cost (USD) \\
\hline
Purity and Abhyankar strategy
  & Claude Fable~5, Claude Code & not recorded & not recorded \\
Example~\ref{ex:cusp}
  & GPT-5.4, ChatGPT Classic & not recorded & not recorded \\
Example~\ref{ex:branch-singular}
  & GPT-6 Astra, Codex desktop, effort xhigh
  & input 133{,}597; cache read 1{,}060{,}352; output 11{,}086
  & 2.95--5.62 \\
Drafting and corrections of this paper
  & Claude Opus~5.5, Claude Code desktop
  & input 282; cache write 553{,}711; cache read 44{,}406{,}202;
    output 130{,}149
  & 15.92 \\
Referee report on the draft
  & Claude Fable~5.1, Claude Code subagent
  & 498{,}610, not split by class & unpriced \\
Comparison of branch loci; Example~\ref{ex:pinched}; its review
  & Claude Opus~5.5, Claude Code desktop and subagent
  & input 304; cache write 888{,}534; cache read 40{,}038{,}008;
    output 214{,}071
  & 19.00 \\
Lean alignment and affine-chart derivation
  & Claude Fable~5.1, Claude Code subagents
  & reported aggregates 265{,}139 and 707{,}872
  & unpriced \\
Actual divisor-and-stalk bridge
  & GPT-5 family
    \cite[\href{https://openai.com/index/introducing-gpt-5/}{model-family page}]{OpenAIGPT5Family25}, Codex desktop
  & final goal telemetry 8{,}663{,}579, unsplit
  & unpriced \\
Paper-one split and mathematical review
  & GPT-6 family, Codex desktop and collaboration agent
  & not exposed & unpriced \\
% >>> [agent-codex 2026-10-28 coauthor-feedback-accounting]
Coauthor-feedback literature and layout revision
  & GPT-6 Sol, Codex desktop, effort xhigh
  & input 178{,}067; cache read 7{,}325{,}312; output 25{,}632
  & API equivalent 2.08 \\
% <<< [agent-codex 2026-10-28 coauthor-feedback-accounting]
\hline
\end{tabular}
\caption{AI usage for this paper.}
\label{tab:tokens}
\end{table}

% >>> [agent-codex 2026-10-28 ai-table-placement]
\FloatBarrier
% <<< [agent-codex 2026-10-28 ai-table-placement]
\bibliographystyle{amsalpha}
\bibliography{pcf}

\newcommand{\etalchar}[1]{$^{#1}$}
\providecommand{\bysame}{\leavevmode\hbox to3em{\hrulefill}\thinspace}
\providecommand{\MR}{\relax\ifhmode\unskip\space\fi MR }
\providecommand{\MRhref}[2]{%
  \href{http://www.ams.org/mathscinet-getitem?mr=#1}{#2}
}
\providecommand{\href}[2]{#2}
\begin{thebibliography}{ASWW26}

\bibitem[{Ant}26a]{Ant26}
{Anthropic}, \emph{{Claude Fable 5}},
  \url{https://www.anthropic.com/news/claude-fable-5-mythos-5}, 2026, Large
  language model, version released 9 June 2026.

\bibitem[{Ant}26b]{AntFable51}
\bysame, \emph{{Claude Fable 5.1}},
  \url{https://www.anthropic.com/claude-fable-and-mythos-5-1}, 2026, Large
  language model, model identifier \texttt{claude-fable-5-1}, used as a Claude
  Code subagent; accessed 4 October 2026.

\bibitem[{Ant}26c]{AntOpus55}
\bysame, \emph{{Claude Opus 5.5}},
  \url{https://www.anthropic.com/claude-opus-5-5}, 2026, Large language model,
  model identifier \texttt{claude-opus-5-5}, used in the Claude Code desktop
  application; accessed 4 October 2026.

\bibitem[Ast20]{Ast20}
Matthieu Astorg, \emph{Dynamics of post-critically finite maps in higher
  dimension}, Ergodic Theory and Dynamical Systems \textbf{40} (2020), no.~2,
  289--308, arXiv:1609.02717.

\bibitem[ASWW26]{ASWW26}
Mohammed Abouzaid, Nikhil Srivastava, Rachel Ward, and Lauren Williams,
  \emph{First proof second batch}, 2026.

\bibitem[BBCL22]{BBCL}
Renee Bell, Jeremy Booher, William~Y. Chen, and Yuan Liu, \emph{Tamely ramified
  covers of the projective line with alternating and symmetric monodromy},
  Algebra Number Theory \textbf{16} (2022), no.~2, 393--446.

\bibitem[BIJ{\etalchar{+}}17]{BIJJLMRS17}
Andrew Bridy, Patrick Ingram, Rafe Jones, Jamie Juul, Alon Levy, Michelle
  Manes, Simon Rubinstein-Salzedo, and Joseph~H. Silverman, \emph{Finite
  ramification for preimage fields of post-critically finite morphisms},
  Mathematical Research Letters \textbf{24} (2017), no.~6, 1633--1647,
  arXiv:1511.00194.

\bibitem[BK10]{BelkKoch10}
James Belk and Sarah Koch, \emph{Iterated monodromy for a two-dimensional map},
  In the Tradition of Ahlfors--Bers, V, Contemporary Mathematics, vol. 510,
  American Mathematical Society, 2010, pp.~1--11.

\bibitem[Col26]{Latexmk26}
John Collins, \emph{Latexmk}, \url{https://ctan.org/pkg/latexmk}, 2026, Version
  4.86a, used with pdfTeX 1.40.27 and TeX Live 2025; accessed 24 September
  2026.

\bibitem[CRS23]{CRS23}
Javier Carvajal-Rojas and Axel St{\"a}bler, \emph{Tame fundamental groups of
  pure pairs and {Abhyankar}'s lemma}, Algebra \& Number Theory \textbf{17}
  (2023), no.~2, 309--358.

\bibitem[DGPS25]{Singular}
Wolfram Decker, Gert-Martin Greuel, Gerhard Pfister, and Hans Sch{\"o}nemann,
  \emph{{Singular} 4-4-1 --- {A} computer algebra system for polynomial
  computations}, \url{https://www.singular.uni-kl.de/}, 2025, Version 4.4.1, as
  distributed with SageMath 10.8; accessed 4 October 2026.

\bibitem[DLTW26]{PcfawdFormalization26}
Taylor Dupuy, Alon Levy, Thomas Tucker, and Zhuoran Wang, \emph{{L}ean
  formalization of {PCF} all the way down},
  \href{https://github.com/tdupu/pcfawd-formalization/tree/05890ba0b0453c289157df45c5bbe472551b1025}{Public
  source repository}, 2026, Registered Palomar source commit
  \texttt{05890ba0b045}.

\bibitem[dMU21]{Lean4}
Leonardo de~Moura and Sebastian Ullrich, \emph{The {L}ean 4 theorem prover and
  programming language}, Automated Deduction --- CADE 28, Lecture Notes in
  Computer Science, vol. 12699, Springer, 2021, pp.~625--635.

\bibitem[DTW26a]{DLTW26Complex}
Taylor Dupuy, Thomas Tucker, and Zhuoran Wang, \emph{{PCF} all the way down
  ai-supplement: The complex analytic case}, AI-generated supplement to
  \emph{{PCF} all the way down}; in preparation, 2026.

\bibitem[DTW26b]{DLTW26Rigid}
\bysame, \emph{{PCF} all the way down ai-supplement: The rigid analytic case},
  AI-generated supplement to \emph{{PCF} all the way down}; in preparation,
  2026.

\bibitem[DTW26c]{DTW26PCFAWD}
\bysame, \emph{{PCF} all the way down: preperiodic subvarieties in any
  codimension, in any characteristic}, Article in preparation, 2026.

\bibitem[Dup26a]{AgentSkillsPublic26}
Taylor Dupuy, \emph{agent-skills-public},
  \url{https://github.com/tdupu/agent-skills-public/tree/1e1b9ba1e6bcf20c31f0bacc799ecdc68b502257},
  2026, Public skill collection at commit
  \texttt{1e1b9ba1e6bcf20c31f0bacc799ecdc68b502257}; checked 28 October 2026.

\bibitem[Dup26b]{Mathcity26}
\bysame, \emph{mathcity: skills and policies for mathematical research with
  {AI} agents},
  \url{https://github.com/tdupu/mathcity/tree/0205e0adfb8963a843c28616fcedceafc7ca5680},
  2026, Commit 0205e0adfb8963a843c28616fcedceafc7ca5680; contains the
  \texttt{using-latexpowers}, \texttt{referee-report},
  \texttt{track-down-reference}, \texttt{update-ai-usage} and
  \texttt{update-tokens} skills.

\bibitem[FS92]{FS92}
John~Erik Forn{\ae}ss and Nessim Sibony, \emph{Critically finite rational maps
  on {$\mathbf{P}^2$}}, The Madison Symposium on Complex Analysis (Madison, WI,
  1991) (Providence, RI), Contemporary Mathematics, vol. 137, American
  Mathematical Society, 1992, pp.~245--260.

\bibitem[GM71]{GM71}
Alexander Grothendieck and Jacob~P. Murre, \emph{The tame fundamental group of
  a formal neighbourhood of a divisor with normal crossings on a scheme},
  Lecture Notes in Mathematics, vol. 208, Springer-Verlag, Berlin--New York,
  1971.

\bibitem[Gro71]{SGA1}
Alexander Grothendieck, \emph{Rev{\^e}tements {\'e}tales et groupe fondamental
  ({SGA} 1)}, Lecture Notes in Mathematics, vol. 224, Springer-Verlag,
  Berlin--New York, 1971.

\bibitem[Hel17]{Hel17}
Paul~Alexander Helminck, \emph{Tropicalizing tame degree three coverings of the
  projective line}, 2017.

\bibitem[HIKT14]{AIM14}
Benjamin Hutz, Patrick Ingram, Sarah Koch, and Thomas Tucker,
  \emph{Postcritically finite maps in complex and arithmetic dynamics}, ARCC
  workshop, American Institute of Mathematics, San Jose, California, 2014,
  Workshop held March 3--7, 2014.
  \url{https://aimath.org/pastworkshops/finitedynamics.html}.

\bibitem[Jon97]{Jon97}
Mattias Jonsson, \emph{Dynamical studies in several complex variables}, Ph.D.
  thesis, Kungliga Tekniska H{\"o}gskolan (KTH), Stockholm, 1997, Paper III is
  the text published in Ergodic Theory and Dynamical Systems 18 (1998)
  171--187; \url{https://www.math.stonybrook.edu/theses/thesis97-1/part1.pdf}.

\bibitem[KS10]{KS10}
Moritz Kerz and Alexander Schmidt, \emph{On different notions of tameness in
  arithmetic geometry}, Mathematische Annalen \textbf{346} (2010), no.~3,
  641--668.

\bibitem[Le23]{Le23}
Van~Tu Le, \emph{Periodic points of weakly post-critically finite all the way
  down maps}, Annali di Matematica Pura ed Applicata \textbf{202} (2023),
  no.~3, 1187--1195, arXiv:2205.03625.

\bibitem[MSP{\etalchar{+}}17]{SymPy17}
Aaron Meurer, Christopher~P. Smith, Mateusz Paprocki, Ondrej Cernik, Sergey~B.
  Kirpichev, Matthew Rocklin, Amit Kumar, Sergiu Ivanov, Jason~K. Moore, Sartaj
  Singh, Thilina Rathnayake, Sean Vig, Brian~E. Granger, Richard~P. Muller,
  Francesco Bonazzi, Harsh Gupta, Shivam Vats, Fredrik Johansson, Fabian
  Pedregosa, Matthew~J. Curry, Andy Terrel, Stepan Roucka, Ashutosh Saboo,
  Isuru Fernando, Sumith Kulal, Robert Cimrman, and Anthony Scopatz,
  \emph{Sympy: symbolic computing in python}, PeerJ Computer Science \textbf{3}
  (2017), e103, Version 1.13.2 used for the analytic example; accessed 24
  September 2026.

\bibitem[{Ope}25]{OpenAIGPT5Family25}
{OpenAI}, \emph{Introducing {GPT-5}},
  \url{https://openai.com/index/introducing-gpt-5/}, 2025, Model-family
  description; the exact model identifier for the divisor-and-stalk session was
  not recorded.

\bibitem[{Ope}26a]{OpenAICodex26}
\bysame, \emph{Codex}, \url{https://openai.com/codex/}, 2026, AI coding agent
  and desktop workflow; accessed 20 September 2026.

\bibitem[{Ope}26b]{OpenAI26}
\bysame, \emph{{GPT-5.4 Thinking}},
  \url{https://openai.com/index/gpt-5-4-thinking-system-card/}, 2026, Large
  language model, version documented 5 March 2026.

\bibitem[{Ope}26c]{OpenAIAstra26}
\bysame, \emph{{GPT-6 Astra}},
  \url{https://developers.openai.com/api/docs/models/gpt-6-astra}, 2026,
  Provider model page for \texttt{gpt-6-astra}.

\bibitem[{Ope}26d]{OpenAISol26}
\bysame, \emph{{GPT-6 Sol}},
  \url{https://developers.openai.com/api/docs/models/gpt-6-sol}, 2026, Provider
  model page for \texttt{gpt-6-sol}.

\bibitem[{Ope}26e]{OpenAIGPT6Family26}
\bysame, \emph{Using {GPT-6}},
  \url{https://developers.openai.com/api/docs/guides/latest-model}, 2026,
  Model-family documentation; the exact deployment identifier for the
  paper-specific Codex session was not exposed.

\bibitem[Rei]{ReidCyclic}
Miles Reid, \emph{Surface cyclic quotient singularities and {Hirzebruch--Jung}
  resolutions}, Author-hosted notes, Example 1.2.

\bibitem[{Sta}18]{Stacks}
The {Stacks Project Authors}, \emph{\textit{Stacks Project}},
  \url{https://stacks.math.columbia.edu}, 2018.

\bibitem[{The}26]{MathlibRev26}
{The mathlib Community}, \emph{mathlib4},
  \url{https://github.com/leanprover-community/mathlib4/tree/56a005855ea17e0b0958f73ec801ae140b0a1679},
  2026, Revision \texttt{56a005855ea1}, tag \texttt{v4.35.0}; the revision used
  for the formalizations; checked 28 October 2026.

\bibitem[TTZ26]{palomar-2026-10-30-000004-v1}
{Taylor Dupuy}, {Thomas Tucker}, and {Zhuoran Wang}, \emph{{PCF all the way
  down}},
  \href{https://palomar-registry.org/entry?id=PALOMAR-2026-10-30-000004&version=1}{Palomar,
  PALOMAR-2026-10-30-000004 v1}, 2026.

\end{thebibliography}

\end{document}